\documentclass[12pt,a4paper]{amsart}

\usepackage{graphicx}
\usepackage{amssymb}
\usepackage{amsmath}

\usepackage{subfig}

\usepackage{booktabs}
\usepackage{multirow}
\usepackage{siunitx}
\usepackage{framed}
\usepackage[pagewise]{lineno}

\usepackage{graphicx}
\usepackage{float,epsfig}
\usepackage[russian,english]{babel}
\usepackage{amsfonts}
\usepackage{amsmath}
\usepackage{subfig}
\usepackage{listings,xcolor}
\graphicspath{{Figures/}}

\newtheorem{theorem}{Theorem}

\newtheorem{defn}{Definition}

\newtheorem{remark}{Remark}
\newtheorem{example}{Example}

\def\Re{\mathop{\rm Re}\nolimits}
\def\Im{\mathop{\rm Im}\nolimits}

\def\Im{\mathop{\rm Im}\nolimits}

\newcommand{\cT}{\mathcal{T}}
\newcommand{\ee}{\varepsilon}
\newcommand{\osc}{\mathrm{osc}}
\newcommand{\symD}{\Omega} 

\newcommand{\symQuad}{Q}
\newcommand{\symQuadC}{\tilde{Q}}

\newcounter{minutes}
\divide\time by 60
\newcounter{hours}
\multiply\time by 60 \addtocounter{minutes}{-\time}
\begin{document}
\def\thefootnote{}

\footnotetext{\texttt{{\tiny File:~\jobname .tex, printed: \number\year-%
\number\month-\number\day, \thehours.\ifnum\theminutes<10{0}\fi\theminutes}}}
\makeatletter

\title[]{CONFORMAL MODULI OF SYMMETRIC \\ CIRCULAR QUADRILATERALS WITH CUSPS}\thanks{The work of the second author was supported by the
Russian Foundation for Basic Research and the Government of the
Republic of Tatarstan, grant No~18-41-160003}

\author[H.~Hakula]{H.~Hakula}
 \address{Aalto University, Espoo, Finland}
 \email[]{Harri.Hakula{@}aalto.fi}

\author[S.~Nasyrov]{S.~Nasyrov}
\address{Kazan Federal University,
         Kazan, Russia}
 \email[]{semen.nasyrov@yandex.ru}

\author[M.~Vuorinen]{M.~Vuorinen}
\address{Department of Mathematics and Statistics, University of Turku,
         Turku, Finland}
\email[]{vuorinen@utu.fi}

\begin{abstract}
We investigate moduli of planar circular quadrilaterals symmetric
with respect to both the coordinate axes. First we develop an analytic
approach which reduces this problem to ODEs and devise a numeric
method to find out the accessory parameters. This method uses the Schwarz
equation to determine conformal mapping of the unit disk onto a
given circular quadrilateral. We also give an example
of a circular quadrilateral for which the value of the conformal modulus can be
found in the analytic form; this example is used to validate the
numeric calculations. We also use another method, so called hpFEM, for the numeric calculation
of the moduli. These two different approaches provide results agreeing
with high accuracy.\end{abstract}

\maketitle

\section{Introduction}\label{intr} \label{1}

A \textit{planar quadrilateral} is a Jordan domain  $Q$ on the
complex plane with four fixed points $z_1$, $z_2$, $z_3$, $z_4$ on
its boundary; we call them vertices of the quadrilateral and
assume that they define positive orientation. If we need to specify the
vertices of a quadrilateral, we write $Q=(Q;z_1,z_2,z_3,z_4)$.  As
well-known, there is a conformal mapping $g$ of $Q$ onto a
rectangle $\Pi=(\Pi;1,1+hi,hi,0)$, $h>0$,  such that the vertices
of $Q$ correspond to the vertices of $\Pi$. Then the value $h$
does not depend on $g$; it is called the conformal modulus of $Q$:
$$
\mbox{\rm Mod}(Q):=h.
$$

Another method, due to L.V. Ahlfors \cite[Thm 4.5, p. 63]{Ah}, to
find the modulus  is to solve  the following Dirichlet-Neumann
boundary value problem for the Laplace equation. Consider a planar
quadrilateral $Q=(Q;z_1,z_2,z_3,z_4)$ with the boundary $\partial
Q = \cup_{k=1}^4\partial Q_k$; all the four boundary arcs are
assumed to be non-degenerate. This problem is
\begin{equation*} 
\left\{\begin{matrix}
    \Delta u& =&\ 0,& \text{on}\ &{\ Q,} \\
    u& =&\ 1,&  \text{on}\ &{\partial Q_1 = (z_1,z_2),}\\
    u& =&\ 0, & \text{on}\ &{\partial Q_3 = (z_3,z_4),}\\
    \partial u/\partial n&  =&\ 0,& \text{on}\ &{\partial Q_2 = (z_2,z_3),}\\
    \partial u/\partial n&  =&\ 0,& \text{on}\ &{\partial Q_4 = (z_4,z_1).}\\
\end{matrix}\right.
\end{equation*}
If we find a solution function $u$ to the above $Q$-problem, then the modulus can be
computed in terms of the solution of this problem as $
\int\!\!\int_Q |\nabla u|^2 dxdy \,.$ We will make use both of the
above two formulations for finding the modulus. The modulus of a quadrilateral is closely
related to the notion of the conformal capacity of a condenser. A condenser in the plane
is a pair $(G,E)$ where $G $ is a domain in the plane and $E$ is its compact subset and 
its capacity is \cite{Du} 
\[  
\inf \int_Q |\nabla u|^2 dxdy \,,   \]
where the infimum is taken over the class of all nonnegative $C^{\infty}(G)$ functions
with  compact support in $G$ and $u(x)\ge 1$ for all $x \in E\,.$

Investigation of conformal moduli of quadrilaterals plays an
important role in  geometric function theory. The method of
conformal moduli is a powerful tool in the theory of
quasiconformal mappings in the plane and in multidimensional
spaces, see \cite{Ah,AB,AVV,Du,hkv,Ku,PoSz}. For instance, many classical
extremal problems of geometric function theory are related to moduli of quadrilaterals
or capacities of condensers \cite{AVV,hkv,Du}.


We note that conformal moduli of quadrilaterals are closely
connected with those of doubly-connected planar domains. Indeed,
all smooth enough doubly-connected domains can be subdivided into
two quadrilaterals. In recent years, a lot of attention has been
paid to numerical computation of conformal moduli of some classes
of quadrilaterals such as those associated with polygonal domains
or domains bounded by circular arcs \cite{BL, cro, gotre1,gotre2,tre1,nas, nasvu}.

We investigate moduli of circular quadrilaterals,
bounded by four circular arcs. Naturally, the vertices of such
quadrilaterals are  the intersection points of the arcs. In
addition, we  will assume that quadrilaterals are symmetric with
respect to the real and imaginary axes and have zero inner angles
at the vertices and that all vertices are on the unit circle. However,
these circular arcs need not be perpendicular to the unit circle.
We also include curvilinear $n$-gons in our examples.

\begin{figure}[h]
    \centering
    \subfloat[{Map of detected fibers with detail area.}]{\includegraphics[width=0.3\textwidth]{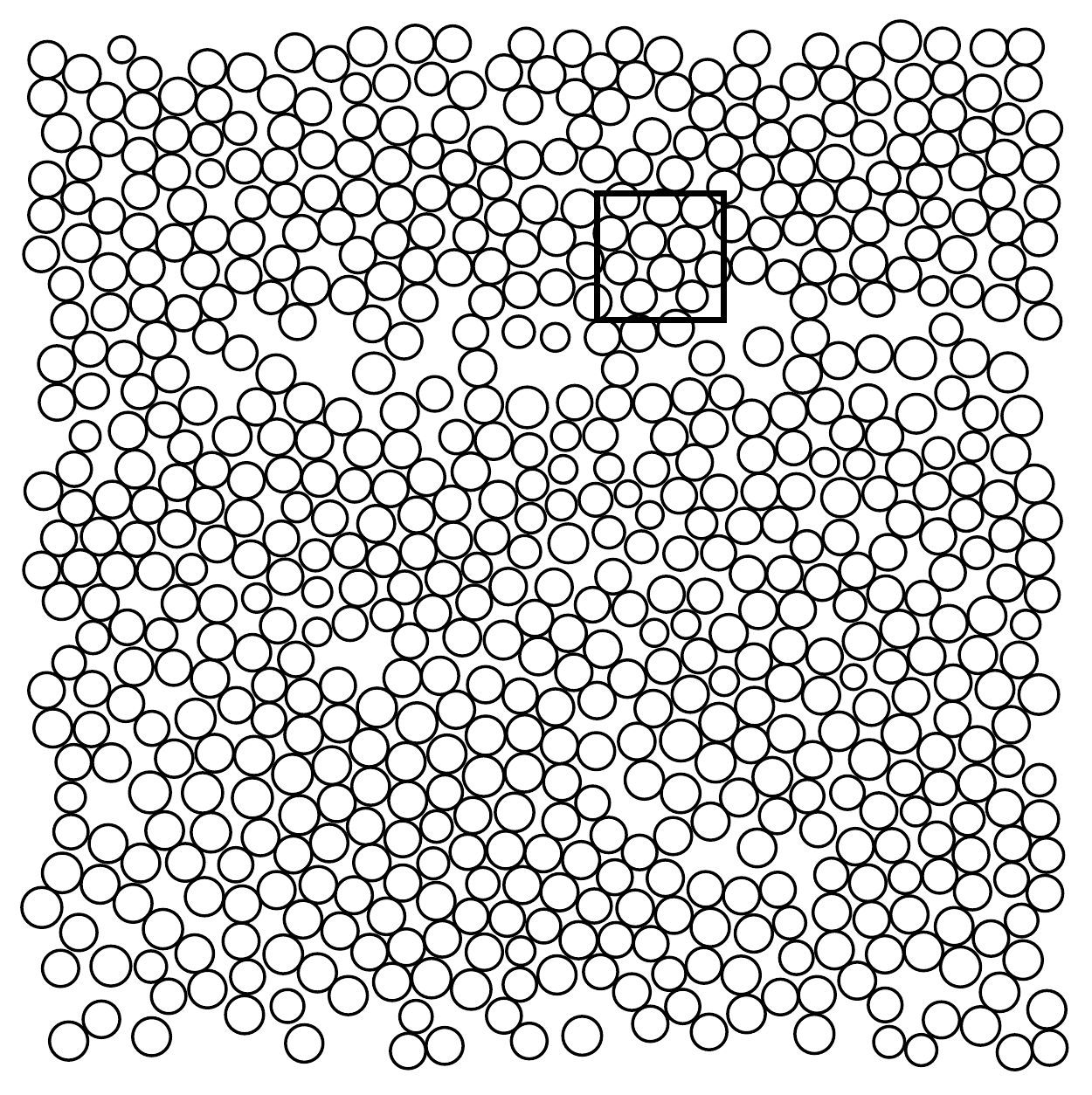}}\quad
    \subfloat[{Detail area with negative colours.}]{\includegraphics[width=0.3\textwidth]{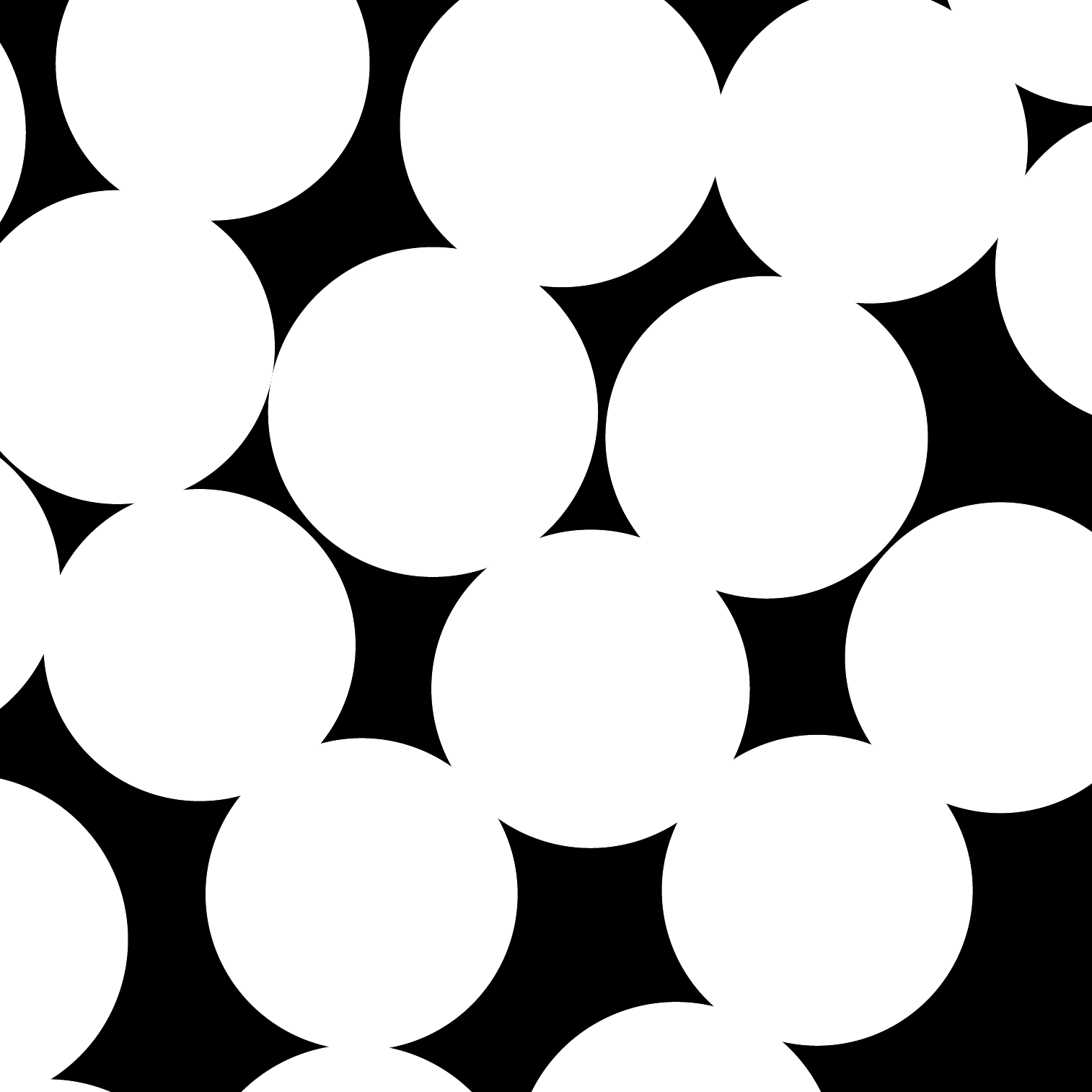}}\quad
    \subfloat[{Detail area after homogenization.}]{\includegraphics[width=0.3\textwidth]{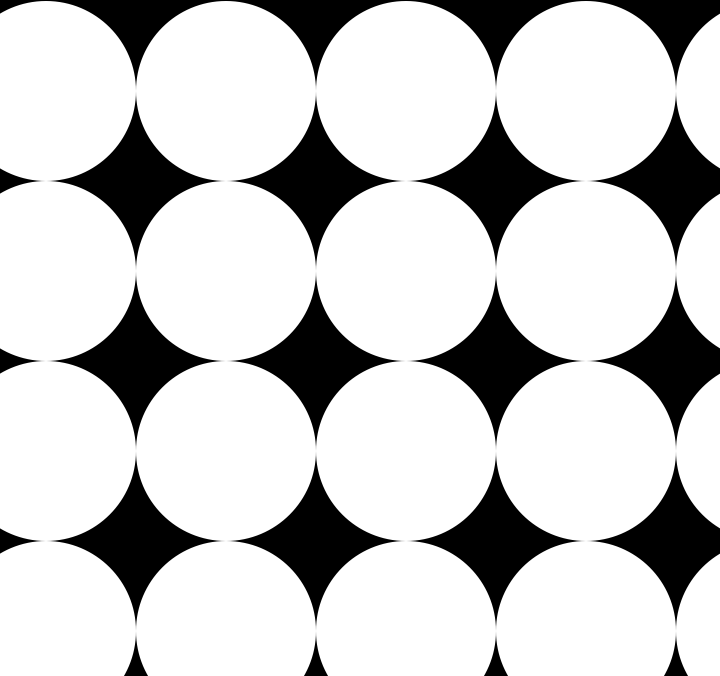}}
    \caption{Circular polygons in applications: Carbon fiber modelling. 
    (a) Map of the measured locations of the fibres within the resin. 
    (b) As the fibres touch, that is, there is contact, the planar intersections of the cavities form circular $n$-gons. 
    (c) After homogenization, the cavities are circular quadrilaterals.
(Data courtesy of I. Babu{\v s}ka, UT Austin.)}\label{fig:appetizer}
\end{figure}
Modelling of carbon fibers induces computational domains that are rich in
such domains \cite{bhl}. In Figure~\ref{fig:appetizer} a map of measured
fiber locations is shown with a detailed image highlighting the domains
bounded by aforementioned circular arcs. Notice that due to measurement tolerances
it would be correct to assume that all sufficiently small gaps could be modelled as closed, i.e.,
neighbouring fibers touching each other. In fact, in order to avoid cusps, in  \cite{bhl}
a minimum distance between fibers was imposed.

Using domain specific discretizations of computational domains as opposed to
traditional triangulations is one of the most active areas of numerical methods
for partial differential equations. In particular, we want to mention
the virtual element method \cite{Beirao2013}
and the cut finite element method \cite{Burman2015}.
In our context, of particular interest is the
contruction of finite elements on curvilinear polygons or $n$-gons \cite{AORW}.
Constructing quadrature rules for such elements is a challenge, and
employing conformal mappings is an intriguing option yet to be fully examined.


Our main goal is to develop numerical methods for calculating
conformal moduli of the kind of circular quadrilaterals and $n$-gons
mentioned above with as high precision as possible. 
Both analytic and purely numerical methods are included in this study.

The analytic method (Section~\ref{2}) uses conformal mappings
of the unit disk onto circular quadrilaterals and their Schwarzian
derivatives. This method is classical;  it is used in many papers
including recent ones. Here we should mention the papers
\cite{brown1,brown2, br_port,kr_port} concerning the usage of
elliptic functions and a spectral Sturm-Liouville problem. The
numerical method reduces to a solution of a pair of ordinary
differential equations (ODE). In Subsection \ref{4}, making use of
the Riemann--Schwarz symmetry principle, we construct  a circular
quadrilateral whose modulus can be determined in the analytic
form. We use this example for testing accuracy of the developed
numerical methods.

The purely numerical method
(Section~\ref{prel}) is based on the $hp$-finite element method (FEM)
implemented by the first author and previously tested in
\cite{hrv,ht}.
In contrast to the first method, the moduli are now computed via potentials
of the associated Dirichlet-Neumann problems.
The $hp$-FEM results are paired with respective
a posteriori error estimates supporting our high confidence in the accuracy
of both methods studied here. Two
error estimators are considered: The $hp$-FEM a posteriori error estimate
based on the auxiliary space methods and 
the physics based reciprocal error \cite{hno,hrv}. Convergence in
the latter, while general, is only a necessary condition and thus
it should always be used in connection with other error
estimators.
The challenges caused by the zero inner angles are well-known. We
deal with this difficulty using geometric mesh grading and control
of the order of the polynomial approximation. 

The two approaches are compared over a parametrized set of 
circular quadrilaterals in the form of graphics and tables. 
The results are in excellent agreement with the
analytic results and support our stated goal of being as accurate as possible.
The $n$-gon test is carried out with the $hp$-version only.
In all cases
exponential convergence is achieved with the $hp$-version at the predicted rates
\cite{schwab}. These observations are supported by both types of
$hp$-error estimators.

We draw our conclusions in Section~\ref{sec:conclusions} and include
a sample implementation of the analytic method in the Appendix.

Several authors have studied various topics about conformal mappings of domains with
circular arc boundaries. The difficulties encountered already in the case of nonsymmetric
quadrilaterals are pointed out in \cite[Section 4.10]{dt}. The interested reader might
want to look at \cite{be,BL, bg,br_port, cro,nas}.

\section{Method of conformal mappings}\label{2}

\subsection{Circular quadrilaterals and the Schwarz
equation}\label{2a}

First we recall some classical results about conformal mapping of
canonical domains onto circular polygons.

Let $D$ be a Jordan domain and let its boundary  consist of $n$
circular arcs $A_{k-1}A_k$, $1\le k\le n$ ($A_0=A_n$). We will
name $D$ a circular polygon; the points $A_k$ are called the
vertices of $D$. Denote by $\alpha_k\pi$ the inner angle of $D$ at
the vertex $A_k$, $0\le\alpha_k\le 2$.

By definition, the Schwarzian derivative of  a meromorphic
function $f$ is the expression
$$
S_f(z)=\left(\frac{f''(z)}{f'(z)}
\right)'-\frac{1}{2}\,\left(\frac{f''(z)}{f'(z)} \right)^2.$$

Let now $f$ be a conformal mapping of the unit disk $U:=\{|z|<1\}$
onto $D$ and denote by $a_k$ the preimage,  under the map $f$, of
$A_k$ lying on the unit circle $\partial U:=\{|z|=1\}$. The following
theorem describes the form of the Schwarzian derivative of~$f$
(see, e.g. \cite[ch.3, \S~1]{gol}, \cite[\S~12]{kop_sht}).

\begin{theorem}\label{sch}
The Schwarzian derivative of the conformal mapping $f$ of $U$ onto
the circular polygon $D$ has the form
\begin{equation}\label{sf}
S_f(z)=\sum_{k=1}^n\frac{(1-\alpha_k^2)/2}{(z-a_k)^2}+\frac{C_k}{z-a_k}\,.
\end{equation}
Here the parameters $C_k$ are some complex numbers satisfying the relations:
$$
\left\{\begin{matrix}
\sum_{k=1}^nC_k=0,\\[2mm]
\frac{1}{2}\sum_{k=1}^n(1-\alpha_k^2)+\sum_{k=1}^nC_ka_k=0,\\[2mm]
\sum_{k=1}^n(1-\alpha_k^2)a_k+\sum_{k=1}^nC_ka_k^2=0.\\
\end{matrix}\right.
$$

\end{theorem}

From Theorem~\ref{sch} we see that the expression (\ref{sf}) for
the Schwarzian derivative of $f$ contains $n$ unknown constants
(or so-called accessory parameters) $C_k$; finding these constants
is a very complicated problem. The problem of determining a
function by its given Schwarzian derivative is well known; many
papers are devoted to this investigation. Various methods are used
to study the problem such as the parametric method
\cite{alex,baib,chist,kolesn,Ku}, boundary value problems
\cite{chibr,tsits1,tsits2,tsits3}, Polubarinova-Kochina's method
\cite{beresl,beres2}, the method of asymptotic integration
\cite{kop_sht}, and others. (Some of these references point out
that the problem of accessory parameters is very important for
investigations in fluid mechanics, especially, in the filtration
theory.)

If we know the values of $C_k$, then the problem of finding $f$ is
reduced to solving the nonlinear third order differential equation
$$
\left(\frac{f''(z)}{f'(z)}
\right)'-\frac{1}{2}\,\left(\frac{f''(z)}{f'(z)}
\right)^2=\sum_{k=1}^n\frac{(1-\alpha_k^2)/2}{(z-a_k)^2}+\frac{C_k}{z-a_k}\,.
$$

The following theorem gives a connection between the problem and
integration of linear second order differential equation (see,
e.g. \cite[ch.VI]{golub}, \cite{kop_sht,tsits3}).

\begin{theorem}\label{sch1}
Let the Schwarzian derivative $S_f$ of $f$ be given. Then $f$ is
defined by $S_f$ up to a M\"obius  transformation. The general solution
of the problem is given by the formula
$$
f(z)=\frac{u(z)}{v(z)}\,.
$$
Here $u$ and $v$ are arbitrary linear  independent solutions of
the equation
\begin{equation}\label{dif2}
h''(z)+(1/2)S_f(z)h(z)=0.
\end{equation}
\end{theorem}


\begin{remark}\label{rem1}
Assume that we seek an odd solution to the problem in a domain
$G$,  containing the origin and symmetric with respect to the
origin, and $S_f(z)$ is an even function in $G$. Then we can take
$u$ and $v$ as odd and even solutions of the equation (\ref{dif2})
in $G$. Therefore, we find $u$ and $v$ as solutions to
(\ref{dif2}) with the following conditions:
\begin{equation}\label{bcond}
u(0)=0,\ \ u'(0)=C\neq 0,\ \ v(0)=1, \ v'(0)=0.
\end{equation}
\end{remark}


\subsection{Conformal mapping of symmetric circular
quadrilaterals}\label{3}

\begin{figure}[ht] \centering
\includegraphics[width=3. in]{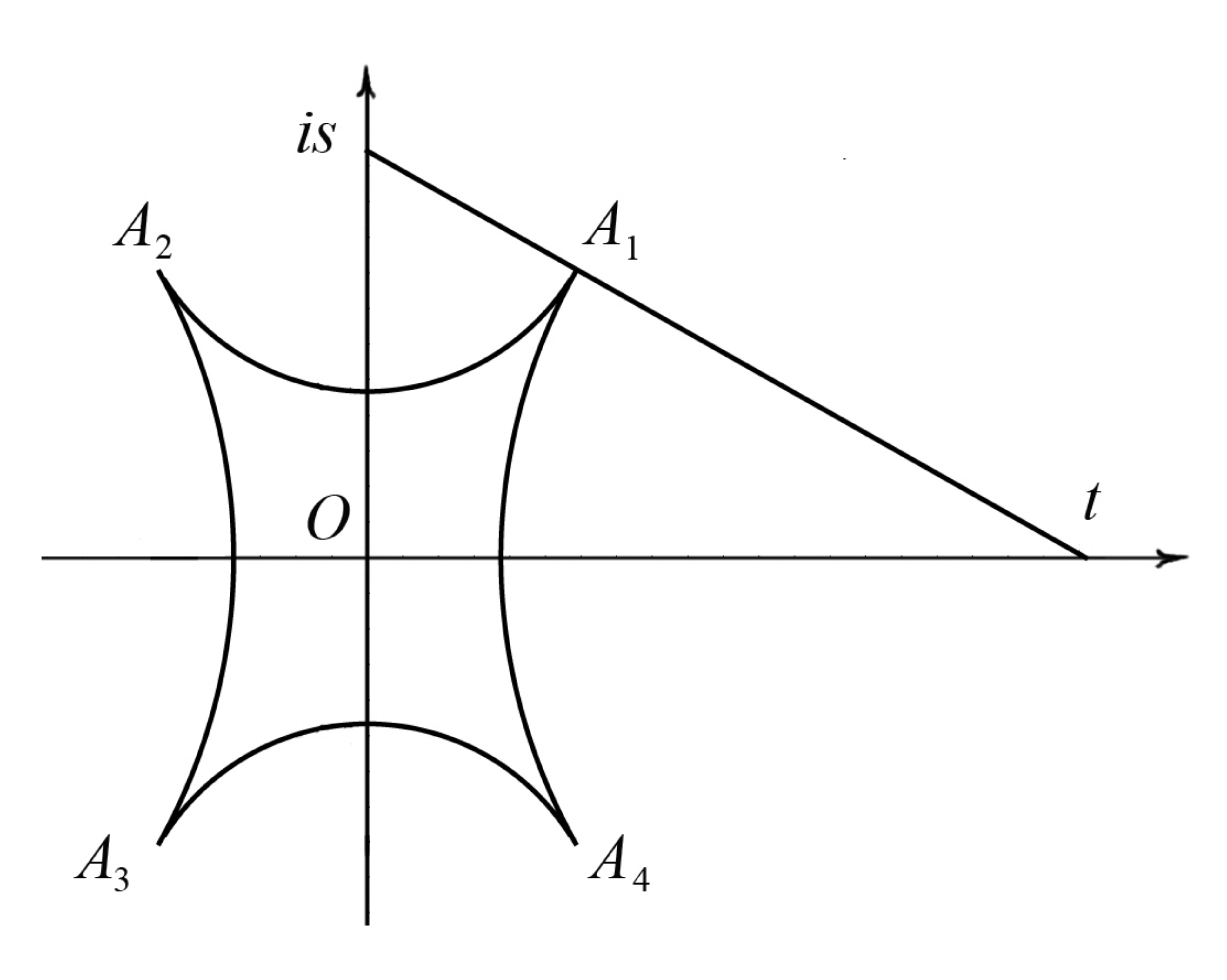}
\caption{Symmetric circular quadrilateral with zero
angles}\label{four}
\end{figure}

We apply Theorems~\ref{sch} and \ref{sch1} in a special case. Let
$A_1A_2A_3A_4$ be a circular quadrilateral with zero inner angles
symmetric with respect to both the axes. Let the centers of the
circles, containing the circular arcs $A_{k-1}A_k$, be at points
$\pm t$, $\pm is$, where $t$, $s>0$. We also assume that at the
points $A_k$ the circles touch each other externally
(Fig.~\ref{four}). Denote the radii of circles centered at $\pm t$
and $\pm is$ by $r_1$ and $r_2$. Then, by Pythagoras' theorem,
$t^2+s^2=(r_1+r_2)^2$.

Denote by $f$ the conformal mapping of the unit disk onto
$A_1A_2A_3A_4$. Because of the symmetry of the quadrilateral with
respect to the coordinate axes and the Riemann--Schwarz symmetry
principle, we may assume without loss of generality that
\begin{equation}\label{sym}
f(\overline{z})=\overline{f(z)}, \quad f(-z)=-f(z).
\end{equation}
Therefore, the preimages $a_k$ of the vertices $A_K$ are symmetric
with respect to the axes, i.e. we can put
\begin{equation}\label{ak}
a_1=e^{i\beta},\ a_2=-e^{-i\beta}, \ a_3=-e^{i\beta}, \  \mbox{\rm
and}\ a_4=e^{-i\beta},\quad 0<\beta<\pi/2.
\end{equation}

Because all the angles of the quadrilateral equal zero, we  have
$\alpha_k=0$, $1\le k \le4$.  Then, by Theorem~\ref{sch}, the
Schwarzian derivative of $f$ has the form
\begin{multline}\label{sdq}
S_f(z)=
\frac{1}{2}\,\left[\frac{1}{(z-e^{i\beta})^2}+\frac{1}{(z+e^{-i\beta})^2}+\frac{1}{(z+e^{i\beta})^2}+\frac{1}{(z-e^{-i\beta})^2}\right]\\
+\left[\frac{C_1}{z-e^{i\beta}}+\frac{C_2}{z+e^{-i\beta}}+\frac{C_3}{z+e^{i\beta}}+\frac{C_4}{z-e^{-i\beta}}\right]\,;
\end{multline}
here the constants $C_k$, $1\le k\le 4$, satisfy
$$
C_1+C_2+C_3+C_4=0,
$$
\begin{equation}\label{cbeta}
e^{i\beta}C_1-e^{-i\beta}C_2-e^{i\beta}C_3+e^{-i\beta}C_4=-2,
\end{equation}
$$
e^{i2\beta}C_1+e^{-i2\beta}C_2+e^{i2\beta}C_3+e^{-i2\beta}C_4=0.
$$
From the first and third equations of the system we obtain
$C_1+C_3=0$, $C_2+C_4=0$, therefore $C_3=-C_1$, $C_4=-C_2$. From
(\ref{sym}) it follows that $C_4=\overline{C_1}$. Moreover,
(\ref{cbeta}) implies $\Re[e^{i\beta}C_1]=-1/2$.

Denote $\Im[e^{i\beta}C_1]=\delta$. Then
$$
e^{i\beta}C_1=-1/2+i\delta,\ e^{-i\beta}C_2=1/2+i\delta,\
e^{i\beta}C_3=1/2-i\delta,\ e^{-i\beta}C_4=-1/2-i\delta.
$$

After simple transformations we obtain (see also \cite{porter})
\begin{equation}\label{swar}
\frac{1}{2}\,S_f(z)=\frac{e^{i2\beta}}{(z^2 -e^{i2\beta})^2}
+\frac{e^{-i2\beta}}{(z^2 -e^{-i2\beta})^2}  - \frac{\gamma} {(z^2
-e^{i2\beta})(z^2 -e^{-i2\beta})}
\end{equation}
where $\gamma=2 \delta\sin 2\beta\in \mathbb{R}$.

By Theorem \ref{sch1}, taking into account (\ref{sym}),  we
represent  $f$ in the form
\begin{equation}\label{fc}
f(z)=C\,\frac{u(z)}{v(z)}\,,\quad C>0,
\end{equation}
where $C$  is a constant and
$$
u(z)=z+\sum_{k=2}^\infty a_k z^k+\ldots\ \mbox{\rm and} \
v(z)=1+\sum_{k=1}^\infty b_k z^k
$$ are linearly independent
solutions to the ODE (\ref{dif2}). Taking into account that, by
(\ref{sym}), $f$ is odd and $S_f$ is an even function, with the
help of Remark~\ref{rem1}, we conclude that $u(z)$ is odd and
$v(z)$ is even. Thus,
\begin{equation}\label{u}
u''(z)+({1}/{2})\,S_f(z)u(z)=0,\quad u(0)=0,\ u'(0)=1,
\end{equation}
\begin{equation}\label{v}
v''(z)+({1}/{2})\,S_f(z)v(z)=0,\quad v(0)=1,\ v'(0)=0.
\end{equation}

Therefore, we have the following result.

\begin{theorem}\label{summ}
Let $f$ be the conformal mapping of the unit disk onto a symmetric
circular quadrilateral $Q$ with zero angles such that the points
$\pm e^ {\pm i\beta}$, $0<\beta< \pi/2$, correspond to the
vertices of $Q$. Then the Schwarzian derivative of $f$ is a
rational function expressed by (\ref{swar}) with some real
$\gamma$, and $f$ has the form (\ref{fc}) where $C$ is a positive
constant,  and the functions $u$ and $v$ are solutions of the
problems (\ref{u}) and (\ref{v}).
\end{theorem}

We should note that the result on the form of the Schwarzian
derivative is actually obtained in \cite[Appendix]{porter} for a
more general case. Here we focus on the case of zero angles.

The equations (\ref{u}) and (\ref{v}) can be used to find the
values of $\beta$ and $\gamma$, corresponding to a given circular
quadrilateral $Q=A_1A_2A_3A_4$.  If we fix some values of
parameters $\beta$ and $\gamma$ and solve the boundary problems
for ODEs, then we find the mapping $f(z)=Cu(z)/v(z)$ up to a
factor $C\neq0$. The obtained function $f(z)=f(z;\beta,\gamma)$
maps the unit disk onto a symmetric circular quadrilateral,
possibly, non-univalently, with zero inner angles. It is evident
that for a given symmetric circular quadrilateral $Q$, there is a
unique pair $(\beta,\gamma)$ such that $f(z;\beta,\gamma)$, with
an appropriate value of $C$, maps the unit disk onto~$Q$.

Therefore, the main problem is to find such a pair $(\beta,\gamma)$
for a given $Q$. We note that the parameter $\beta$ has a very
simple geometric meaning. Finding $\beta$ is equivalent to  finding
the conformal modulus of $Q$. Actually, because of the property of
conformal invariance, the modulus of $Q$ is equal to the modulus
of the unit disk with vertices (\ref{ak}) which depend only on
$\beta$. The parameter $\gamma$ has no simple geometric meaning but
it also affects the geometry of $Q$.


To find $(\beta,\gamma)$ numerically, we seek
$u(e^{i\theta})$ and $v(e^{i\theta})$ as solutions of the following boundary
value problems for ODEs:
\begin{equation}\label{ut}
u''(re^{i\theta})+\frac{e^{i2\theta}}{2}\,S_f(re^{i\theta})u(re^{i\theta})=0,\
0\le r\le 1,\quad u(0)=0,\ u'(0)=e^{i\theta},
\end{equation}
\begin{equation}\label{vt}
v''(re^{i\theta})+\frac{e^{i2\theta}}{2}\,S_f(re^{i\theta})v(re^{i\theta})=0,\
0\le r\le 1,\quad v(0)=1,\ v'(0)=0,
\end{equation}
and determine $f(e^{i\theta})=C u(e^{i\theta})/v(e^{i\theta})$.
Then we determine the values $T=f(1)$, $S=f(i)/i$, and the values
of the radii $R_1$  and $R_2$ of the circles containing the circular
arcs $f(e^{i\theta})$, $-\beta<\theta<\beta$, and
$f(e^{i\theta})$, $\beta<\theta<\pi-\beta$. All these values
depend on the parameters $\beta$ and $\gamma$. Then we compare the
ratios $S/T$ and $R_2/R_1$ with the given ones, $s/t$ and
$r_2/r_1$. Therefore, we have two equations to determine $\beta$
and $\gamma$:
\begin{equation}\label{ST}
\frac{S(\beta,\gamma)}{T(\beta,\gamma)}=\frac{s}{t}\,,\quad
\frac{R_2(\beta,\gamma)}{R_1(\beta,\gamma)}=\frac{r_2}{r_1}\,.
\end{equation}
This system has a unique solution and the obtained value of
$\beta$ enables us to find the modulus of~$Q$.

Because our  main  goal is to determine the conformal modulus of
$Q$ solving the system (\ref{ST}) and because the ratios $R/T$ and
$R_2/R_1$ do not change under homotheties, below we will often
assume that the constant  $C$ in (\ref{fc}) equals $1$.

\subsection{Example of circular quadrilateral with exactly known
modulus}\label{4}

Unfortunately, there are very few examples of concrete circular
quadrilaterals with exactly known values of the conformal modulus.
In this section, with the help of the Riemann--Schwarz symmetry
principle and the Schwarz--Christoffel formula, we give an example
of this type.

\begin{figure}[ht] \centering
\includegraphics[width=5.0 in]{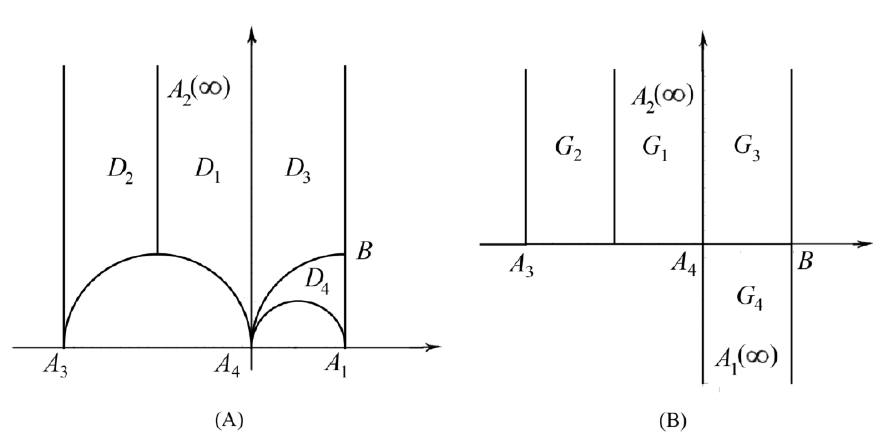}\
\caption{(A) Circular quadrilateral in the upper half-plane; \ (B)
its strip-shaped image.}\label{four1}
\end{figure}

Consider the circular quadrilateral $Q$, $Q\subset U$,  with
vertices lying on the unit circle $\partial U$ at the points
$$
A_1=e^{i\alpha},\ A_2=-e^{-i\alpha}, \ A_3=-e^{i\alpha}, \
\mbox{\rm and}\ \ A_4=e^{-i\alpha},\quad
\alpha=\arcsin{(1/\sqrt{3})}.
$$
Let the boundary arcs of $Q$ be orthogonal to the unit circle.
The M\"obius  transformation
$$
\omega=(i\cot \alpha)\,\,\frac{e^{-i\alpha}-z}{e^{-i\alpha}+z}\,
$$
maps $Q$ conformally onto the circular quadrilateral $D$ lying in
the upper half-plane of the variable $\omega$, bounded by two
rays, $\{\Re \omega=-2, \Im \omega\ge 0\}$, $\{\Re \omega=1, \Im
\omega\ge 0\}$, and two semicircles, $\{|\omega+1|=1, \Im
\omega\ge 0\}$ and $\{|\omega-1/2|=1/2, \Im \omega\ge 0\}$
(Fig.~\ref{four1}~(A)). For convenience and brevity of notation,
we use the same notations for boundary points corresponding to
each other in different complex planes under the applied conformal
mappings.

Denote by $D_1$ the subdomain of $D$ lying in the strip $\{-1<\Re
\omega<0\}$.

Let $F$ be the conformal map of $D_1$ onto the half-strip
$G_1=\{-1<\varphi<0,\psi>0\}$ in the $w$-plane ($w=\varphi+i
\psi$) such that $F(-1+i)=-1$, $F(0)=0$, $F(\infty)=\infty$.

Applying the Riemann-Schwarz symmetry principle, we extend  the mapping $F$ to the
domains $D_2$ and $D_3$, symmetric to $D_1$ with respect to lines
$\Re \omega=-1$  and $\Re \omega=0$, resp. Then the extended mapping,
for which we keep the same notation $F$, maps the union of
domains $D_1\cup D_2\cup D_3$ (supplemented with their common
boundary arcs) onto the strip $\{-2<\varphi<1\}$ consisting of
three half-strips, $G_1$, $G_2$, and $G_3$ (Fig.~\ref{four1}). At
last, we can extend, by symmetry,  $F$ to the domain $D_4$,
symmetric to $D_3$ with respect to the boundary arc $A_4B$, lying
on the unit circle $\{|w|=1\}$. The extended function maps
conformally $D_4$ onto the half-strip $G_4$ symmetric to $G_4$
with respect to the real axis. As a result, we conclude that the
domain $D$ is conformally equivalent to the strip-shaped domain
$G$ glued from the half-strips $G_k$, $1\le k\le 4$, along their
common boundary segments.

Let us map conformally the upper half-plane in the $\zeta$-plane
onto $G$ such that the points $-1/\lambda$, $-1$, $1$, and
$1/\lambda$ ($\lambda>1$) correspond to $A_2$, $A_3$, $A_4$, and
$A_1$. The desired mapping is given by the Schwarz-Christoffel
integral
$$
G(\zeta)=c\int_{1}^{\,\zeta}
\sqrt{\frac{1-t}{1+t}}\,\,\frac{dt}{1-\lambda^2 t^2}
$$
with some constant $c>0$. In a neighborhood of $\zeta=1/\lambda$
we have
\begin{equation}\label{as1}
G(\zeta)\sim i\frac{c}{2\lambda} \sqrt{\frac{1-\lambda}{1+\lambda
}}\,\log(\zeta-1/\lambda).
\end{equation}
By a similar way, as $\zeta\to -1/\lambda$, we have
\begin{equation}\label{as2}
G(\zeta)\sim -i\frac{c}{2\lambda} \sqrt{\frac{1+\lambda}{1-\lambda
}}\,\log(\zeta+1/\lambda).
\end{equation}

We recall that the function $\log(\zeta)$,  maps conformally the
upper half of a sufficiently small neighborhood of the origin onto
a half-strip-like domain of width $\pi$.

Taking into account the values of widths of the half-strip parts
of $G$ and the asymptotics (\ref{as1}), (\ref{as2}), we obtain
$$
\frac{\pi c}{2\lambda} \sqrt{\frac{1-\lambda}{1+\lambda }}\,=1,
\quad   \frac{\pi c}{2\lambda} \sqrt{\frac{1+\lambda}{1-\lambda
}}\,=3,
$$
therefore, $(1+\lambda)/(1-\lambda)=3$ and $\lambda=1/2$.

Now we can find the value of the conformal modulus of $Q$. Because
of the invariance of the modulus under conformal mappings, we see
that $\mbox{\rm Mod}(Q)$ is equal to the modulus of  the quadrilateral
which is the upper half-plane with vertices $\pm 1$, $\pm
1/\lambda$. Therefore, it can be computed via elliptic integrals
(see, e.g. \cite{akhiezer}, \cite{AVV}):
$$
(\mbox{\rm
Mod}(Q))^{-1}=\frac{2K(\lambda)}{K(\lambda')}
$$
where
$$
K(\lambda)=\int_0^1\frac{dt}{\sqrt{(1-t^2)(1-\lambda^2t^2)}}
$$ is the complete elliptic integral of the
first kind and $\lambda'=\sqrt{1-\lambda^2}=\sqrt{3}/2$. At last,
we obtain
$$
\mbox{\rm
Mod}(Q)=\frac{K(\sqrt{3}/2)}{2K(1/2)}\,=0.6396307855855...
$$

Now we will find  the Schwarzian derivative of the conformal
mapping of the unit disk onto $Q$. As we noted in
Subsection~\ref{3}, it has the form (\ref{swar}). The cross-ratio
of a quadruple,
$$(z_1,z_2;z_3,z_4)=\frac{(z_3-z_1)(z_4-z_2)}{(z_3-z_2)(z_4-z_1)}$$
is invariant under M\"obius  transformations. Comparing the
cross-ratios for vertices of the two quadrilaterals, the first one
of which is the unit disk with vertices $\pm e^{\pm i\beta}$ and
the second one is the upper half-plane with vertices $\pm 1$, $\pm
1/\lambda$ we obtain
$$\sin\beta=\frac{1-\lambda}{1+\lambda}\, =\frac{1}{3}\,,$$
therefore, $\beta=\arcsin(1/3)$.

Careful analysis of the Schwarzian derivative in small
neighborhoods of the points $\pm e^{\pm i\beta}$ shows that, in
the considered case, the parameter $\gamma$ in (\ref{swar}) equals
$2/3$.

\subsection{Conformal mapping of circular $n$-gons}

It is of interest to consider circular $n$-gons with $n>4$ and
compute moduli of quadrilaterals which are obtained from them
after fixing four of their vertices.

Here we give some examples of circular $n$-gons with zero angles
and known conformal moduli of quadrilaterals constructed on the
base of these $n$-gons; they can be also used for testing the
error of the $hp$-FEM in finding conformal moduli
(Section~\ref{prel}).\medskip

\begin{example}\label{exam1}
In the $z$-plane, $z=x+iy$, we consider a circular hexagon $H$
with zero angles.  The hexagon is obtained from the half-strip
$\{-2<x<2,y>0\}$ by removing points lying in the disks
$\{(x+1)^2+y^2\le 1$,  $\{(x-1/3)^2+y^2\le 1/9\}$,
$(x-5/6)^2+y^2\le 1/36\}$, $\{(x-3/2)^2+y^2\le 1/4\}$. It has
vertices at the points (Fig.~\ref{fig:ngon}~(A))
\[A(-2,0),\, B(0,0), \, C(2/3,0), \,D(1,0),\,
E(2,0), \,\,{\rm and}\,\, F(\infty)\,.\]

\begin{figure}[ht] \centering
\includegraphics[width= 5. in]{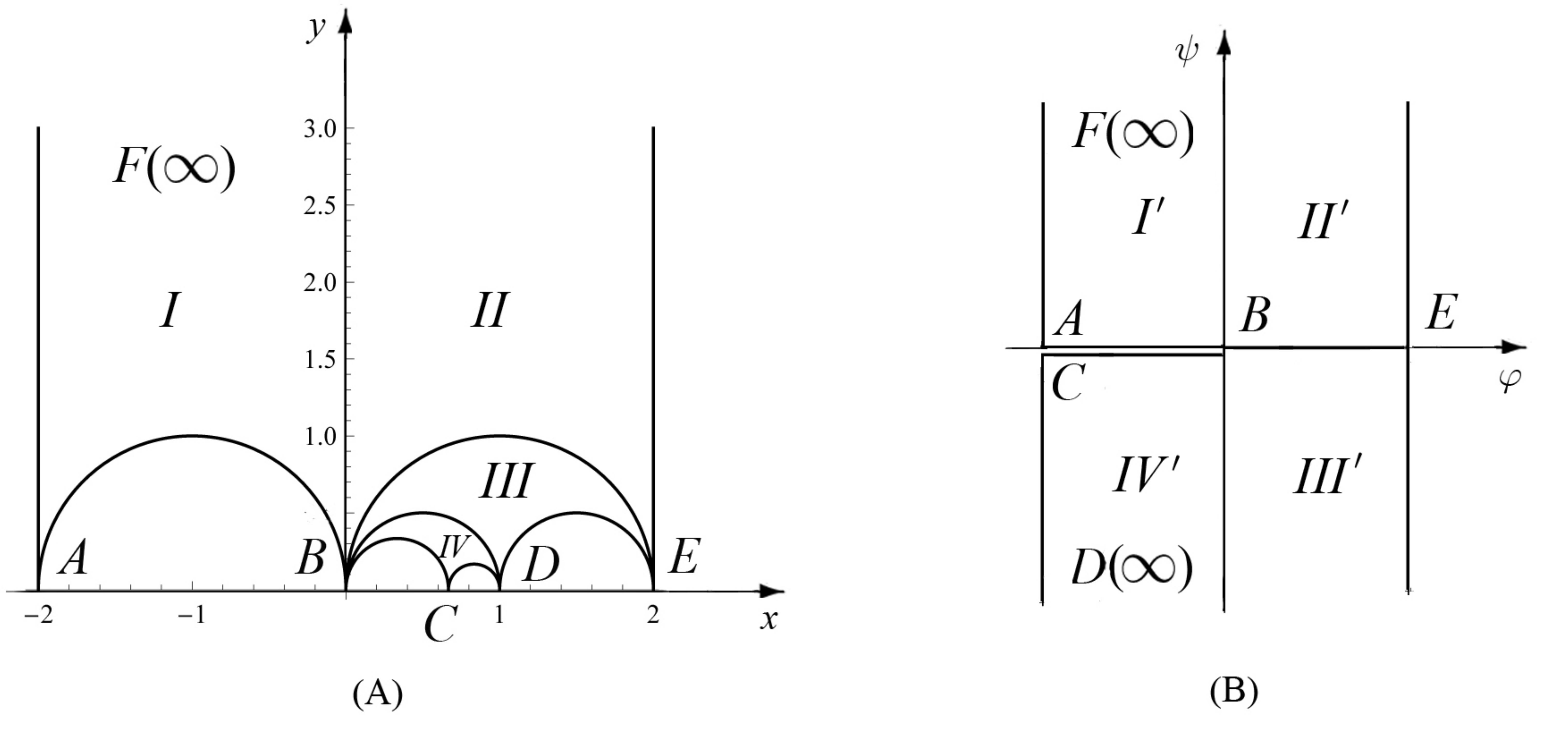}
\caption{(A) Circular hexagon in the upper half-plane; \ (B) its
polygonal image.}\label{fig:ngon}
\end{figure}
\medskip

Let us map conformally the circular triangle
$I:=\{-2<x<0,y>0\}\setminus \{(x+1)^2+y^2\le 1\}$ onto the
half-strip $I':=\{-\pi/2<\varphi<0,\psi>0\}$ in the $w$-plane,
$w=\varphi+i \psi$. Applying three times the Riemann-Schwarz
symmetry principle, we extend the mapping step by step to the
circular triangles designated on the Fig.~\ref{fig:ngon}~(A), by
$I\!I$, $I\!I\!I$, and $IV$. The extended function $f$ maps the
triangles onto the half-strips $I\!I'$, $I\!I\!I'$, and $IV'$
(Fig.~\ref{fig:ngon}~(B)).

Therefore, $f$ maps conformally $H$ onto the strip
$\{|\varphi|<\pi/2\}$ with the slit along the segment
$\{-\pi/2<\varphi<0,\, \psi=0\}$. Denote this domain by $\Omega$.
The function $\omega=g(w):=2\sin w-1$ maps conformally $\Omega$
onto the $\omega$-plane, $\omega=\xi+ i \eta$ with two slits along
the rays $\{\xi\le -1,\, \eta=0\}$ and $\{\xi\ge 1,\, \eta=0\}$.
Then we apply the function inverse to the Joukowsky function:
$\sigma=h(\omega)=\omega+\sqrt{\omega^2-1}$ with an appropriate
choice of regular branch of the square root. The composition
$h\circ g\circ f$ maps conformally the hexagon $H$ onto the upper
half-plane with the following correspondence of the points:
\begin{equation}\label{hexag}
A\mapsto -(3+2\sqrt{2}), B\mapsto -1, C\mapsto -(3-2\sqrt{2}),
D\mapsto 0, E\mapsto 1, F\mapsto \infty.
\end{equation}
Therefore, we have the following result.\medskip

\textit{The hexagon $H$ is conformally equivalent to the upper
half-plane with the correspondence of points given by
(\ref{hexag}).}\medskip

Because the conformal modulus of the quadrilateral which is the
upper half-plane with four fixed vertices on the real axis is
well-known, we can fix any four of the six vertices and easily
compute the modulus of the obtained quadrilateral.
\end{example}

\begin{remark}\label{rem2}
In Section~\ref{prel} we use this example to verify the accuracy
of the $hp$-FEM for determining conformal moduli of circular
$n$-gons. This method needs calculation of double integrals over a
given $n$-gon, therefore, it is better to apply it in a bounded
domain. Since the modulus is a conformal invariant, we can
consider, instead of the unbounded hexagon $H$, its conformal
image under the M\"obius transformation
$$
w=T(z)=\frac{4-(1-3i)z}{4-(1+3i)z}\,.
$$
This transformation maps the upper half-plane onto the unit disk;
and $Q$ corresponds to the hexagon bounded by circular arcs
orthogonal to the unit circle. Moreover, there is the following
correspondence between the vertices of $Q$ and their images:
\begin{equation*}
A\mapsto -i,\ B\mapsto 1,\ C\mapsto \frac{8+15i}{17}\,,\ D\mapsto
i,\ E\mapsto \frac{-4+3i}{5}\,,\ F\mapsto \frac{-4-3i}{5}\,.
\end{equation*}
The obtained hexagon is considered in Subsubsection~\ref{exam}
(see Fig.~\ref{fig:CNG}).
\end{remark}

\medskip

\begin{example}\label{exam2}
Given $n\ge 4$, consider the circular $n$-gon $P_n$ which is
obtained from the half-strip
$$\{0<x<2(n-2),y>0\}$$ by removing the disks
$D_k:=\{(x-(2k-1))^2+y^2\le 1\}$, $1\le k\le n-2$. It has zero
angles and vertices at the points $0$, $2$, $4$,\ldots,$2n-2$,
$\infty$. We map the triangle $\{0<x<2,y>0\}\setminus D_1$ onto
the half-strip $\{0<\varphi<2, \psi>0\}$ and extend the mapping by
symmetry to a conformal mapping of the $n$-gon onto the half-strip
$S:=\{0<\varphi<2(n-2), \psi>0\}$. Then we map $S$ onto the upper
half-plane by the function $\zeta=-\cos \frac{\pi w}{2(n-2)}$.
Then the vertices of $P_n$ 
are mapped to the points $-\cos \frac{\pi k}{n-2}$,
$k=0,1,\ldots,n-2$, and $\infty$. As in Example~1, fixing four
vertices of $P_n$ we can find exact value of modulus of the
obtained quadrilateral.
\end{example}

\subsection{Numeric results}\label{5}


Here we  give a numerical algorithm  to find the values of $\beta$
and $\gamma$ in \eqref{swar} for a given symmetric quadrilateral
with zero inner angles (see Fig.\ref{four}). For this, as
indicated above, we need to solve the system~\eqref{ST}.

Denote $k:=s/t$, $K:=r_2/r_1$. Then \eqref{ST} has the form
\begin{equation}\label{sys2}
S(\beta,\gamma)/T(\beta,\gamma)=k,\quad
R_2(\beta,\gamma)/R_1(\beta,\gamma)=K.
\end{equation}

First we describe how, for given arbitrary $\beta$ and $\gamma$,
to find the centers, $T(\beta,\gamma)$ and $iS(\beta,\gamma)$, and
the radii, $R_1(\beta,\gamma)$ and $R_2(\beta,\gamma)$, of circles
which contains the boundary circular arcs of the corresponding
circular polygon. 
We note that boundary arcs are symmetric with respect to either
the real or the imaginary axis, therefore, it is sufficient to
determine two distinct points for each of the circles. We can find
$f(e^{i\theta})=u(e^{i\theta})/v(e^{i\theta})$, solving the
equations (\ref{ut}) and (\ref{vt}), where $S_f$ is defined by
\eqref{swar} and corresponds to the fixed values of $\beta$ and
$\gamma$. Then we take two different values of $\theta$ from $[0,
\beta)$, say, $\theta_1=0$ and $\theta_2=\beta/2$, and find the
values of $f(e^{i\theta_1})$ and $f(e^{\theta_2})$. Let
$x_1=f(e^{i\theta_1})$ and $x_2+iy_2=f(e^{\theta_2})$. Next we
determine
$$
T = (1/2) (x_1 + x_2 + y_2^2/(x_2 - x_1)),\ R_1=|t-x_1|
$$ where
$x_1=f(e^{i\theta_1})$ and $x_2+iy_2=f(e^{\theta_2})$. By a
similar way, we fix two angles in $(\beta,\pi/2]$, say,
$\theta_3=\pi/4+\beta/2$ and $\theta_4=\pi/2$, and find
$$
S = (1/2) (y_3 + y_4 + x_4^2/(y_4 - y_3)),\ R_2=|s-y_4|,
$$
where $x_3+iy_3=f(e^{i\theta_3})$, $iy_4=f(e^{\theta_4})$.

Now we describe how to solve the system \eqref{sys2}. Initially,
for a given fixed $\beta$, we solve the first equation from
(\ref{sys2}) with respect to  $\gamma$. We should note that for
given values of $\beta$ and $k$, a symmetric circular
quadrilateral and, therefore, $\gamma$ are not uniquely
determined. We consider circular polygons such that the circles,
containing their boundary, touch each other externally at
intersection points. But there could be another circular polygon
with the same values of $\beta$ and $k$ and the circles touching
internally. To avoid this, first, for a given $\beta$ we need to
determine the values of $\gamma$, $A_\gamma$ and $B_\gamma$,
$A_\gamma<B_\gamma$, for which we obtain circular quadrilaterals
with a pair of sides lying on parallel straight lines. The values
$A_\gamma$ and $B_\gamma$ correspond to the conditions $x_1=0$
($R_1=\infty$) and $y_4=0$ ($R_2=\infty$).

We find the values of $A_\gamma$ and $B_\gamma$ by the bisection
method on some segment $I=[a_\gamma, b_\gamma]$. The segment $I$
must be sufficiently large and contain $A_\gamma$ and
$B_\gamma$. Using a numerical experiment, we determined that for a
wide class of $\beta$ and $k$ the following values of the
parameters are appropriate:
$$
a_\gamma=0.7 - (4/\pi) \beta, \quad  b_\gamma = 1.2 -
(3/\pi)\beta.
$$
When the values of $A_\gamma$ and $B_\gamma$ are found, we
determine the desired value of $\gamma=\gamma(\beta)$.

To fulfill the second equality in (\ref{sys2}), we solve the
equation
$$
R_2(\beta,\gamma(\beta))/R_1(\beta,\gamma(\beta))=K
$$
making use of the bisection method on the segment $[0,\pi/4]$. We
do not know a priori, whether $\beta<\pi/4$ (this means that the
modulus is less than $1$). If it turns out that the numerical
value of the desired modulus is greater than $1$ and, therefore,
$\beta\ge\pi/4$, then the bisection method on the segment
$[0,\beta/4]$ converges to the boundary value $\beta=\pi/4$. In
this case, we swap the values of $t$ and $s$, as well as $r_1$ and
$r_2$, and repeat the calculations for these updated values; at
the end, the found value of $\beta$ must be changed to
$\pi/2-\beta$.

For numeric calculations we used the Wolfram Mathematica software.
If we want to obtain the approximate value of the modulus quickly
and with accuracy about $10^{-6}$, for solving Cauchy's problems
for differential equations with the help of NDSolve we can use the
option 'PrecisionGoal$-\!\!\!>$15. To find the parameters
$A_\gamma$ and $B_\gamma$ with the help of the bisection
method, it is sufficient to use 10 iterations; and for each of the
parameters, $\beta$ and $\gamma$, we used 25 iterations.

In Appendix~A we give the Mathematica code for calculation of
conformal moduli. The input values of $t$, $s$,  $r_1$, and $r_2$
(lines 1--4) match to the example 2) below with $\alpha=\pi/5$,
$j=3$. The output is the found values of $\mbox{\rm Mod}(Q)$,
$\beta$ and $\gamma$ 
(line 118).

If we need a higher accuracy, we can first find the approximate
values $\beta$ and $\gamma$ with accuracy $10^{-6}$. Denote them
by $\beta_0$ and $\gamma_0$. Then we use the bisection method with
respect to $\beta$ and $\gamma$ assuming that
$\beta\in[\beta_0-\varepsilon,\beta_0+\varepsilon]$ and
$\gamma\in[\gamma_0-\varepsilon,\gamma_0+\varepsilon]$ with
sufficiently small $\varepsilon$, say $\varepsilon=2\cdot10^{-6}$.
Certainly, in the case, we omit the first two steps connected with
finding the values of $A_\gamma$ and $B_\gamma$. To find $\beta$
and $\gamma$, we use NDSolve with option
'PrecisionGoal$-\!\!\!>$30; for each of the parameters, the number
of iterations is $30$. With this enhanced method, the accuracy is
about $10^{-10}-10^{-11}$ though it needs much more computing time
(a few minutes instead of 10--15 seconds).
\medskip

Now we give numerical results.\medskip

1) For $k=\sqrt{2}$ and $K=2$ we know the exact values
$\beta=\arcsin(1/3)$, $\gamma=2/3$ (see Subsection~\ref{4}).

Using the options PrecisionGoal$-\!\!\!>$15,
WorkingPrecision$-\!\!\!>$30 with the number of steps  equals 30
for each of the parameters, we find that  the approximate values
are
$$\sin\beta=0.333 333 333 333 244 1,\quad \gamma=0.666 666 666 666
788.$$ Therefore, the absolute error is about $1.2\cdot 10^{-13}$.
\medskip

2) Assume that the vertex $A_1$ is $e^{i\alpha}$,
$\alpha=\pi/n$, $n\in \mathbb{N}$, $4\le n\le 8$
(Fig.\ref{four}). For every $n$ we consider  the following five
values of $t$:
$$
t=1+0.2 j (1/\cos \alpha-1),\quad 1\le j\le 5.
$$
Then
$$
s=\frac{t\sin\alpha}{t-\cos\alpha}\quad r_1=|e^{i \alpha}-t|,
\quad r_2=|e^{i \alpha}-is|.
$$

We computed the values of conformal modulus for these 25 cases. In
Table~\ref{tab2} 
we give the values obtained with rough
accuracy, with higher accuracy, and by the $hp$-FEM (see
Section~\ref{prel}). We see that for given $\theta$, the
difference in results, given in the fourth (higher accuracy) and
fifth ($hp$-FEM variant) columns does not exceed $5\cdot
10^{-10}$. This indicates fairly good accuracy of the suggested
methods.

We note that in 
Subsubsection~\ref{exam} two the following cases are considered in
more detail: $n=4$, $j=1$ (quadrilateral $Q_1$) and $n=8$, $j=5$
(quadrilateral $Q_2$).

\begin{center}
\begin{table}[ht]
\caption{The values of moduli of circular quadrilaterals.} {\small
\centering
\begin{tabular}{|c|c|c|c|c|}
  \hline
  $\alpha$&j&rough  accuracy&higher accuracy&$hp$-FEM \\
  \hline
\multirow{5}{*}{$\pi/4$}
&1&$\phantom{h}1.65195637087856\phantom{h}$& $\phantom{h}1.65195641811156\phantom{h}$&$\phantom{h}1.65195641811801\phantom{h}$\\
  \cline{2-5}
&2&$1.41312892318176$& $1.41312882432748$& $1.41312882433334$\\
  \cline{2-5}
&3&$1.23851630005081$ & $1.23851628549016$& $1.23851628549600$\\
  \cline{2-5}
&4&$1.10517568205164$& $1.10517573064876$& $1.10517573065505$\\
  \cline{2-5}
&5&$\!\!\!\!\!\!1.$ (sharp value)& $\!\!\!\!\!\!1.$ (sharp value)&$1.00000000000704$\\
  \hline
\multirow{5}{*}{$\pi/5$}
&1&$0.98160716203795$& $0.98160730939538$& $0.98160730941547$\\
  \cline{2-5}
&2&$0.88131392216865$& $0.88131392866493$& $0.88131392869094$\\
  \cline{2-5}
&3&$0.79679231514866$& $0.79679236427334$& $0.79679236430546$\\
  \cline{2-5}
&4&$0.72458905475484$& $0.72458889240001$& $0.72458889243949$\\
  \cline{2-5}
&5&$0.66218846198336$& $0.66218813398119$& $0.66218813402464$\\
  \hline
\multirow{5}{*}{$\pi/6$}
&1&$0.69813401618400$&$0.69813355689778$&$0.69813355697485$\\
  \cline{2-5}
&2&$0.63911291428315$&$0.63911229266297$&$0.63911229274088$\\
  \cline{2-5}
&3&$0.58614443760266$&$0.58614411420414$&$0.58614411428162$\\
  \cline{2-5}
&4&$0.53833141064728$&$0.53833144748697$&$0.53833144756331$\\
  \cline{2-5}
&5&$0.49493995006987$&$0.49493951440663$& $0.49493951447948$\\
  \hline
\multirow{5}{*}{$\pi/7$}
&1&$0.54204363753707$&$0.54204377899126$&$0.54204377906567$\\
  \cline{2-5}
&2&$0.50133118737030$&$0.50133063755764$&$0.50133063763325$\\
  \cline{2-5}
&3&$0.46350927171723$&$0.46350872114770$&$0.46350872122462$\\
  \cline{2-5}
&4&$0.42826417448376$&$0.42826373909062$&$0.42826373916846$\\
  \cline{2-5}
&5&$0.39531876405162$&$0.39531863465020$&$0.39531863472915$\\
   \hline
\multirow{5}{*}{$\pi/8$}
&1&$0.44327621319647$& $0.44327582367411$&$0.44327582393810$\\
  \cline{2-5}
&2&$0.41254694695236$&$0.41254658974644$&$0.41254659003158$\\
  \cline{2-5}
&3&$0.38338345187308$&$0.38338339855016$&$0.38338339885322$\\
  \cline{2-5}
&4&$0.35565053319540$&$0.35565066792949$&$0.35565066823961$\\
  \cline{2-5}
&5&$0.32922105387009$&$0.32922144646084$&$0.32922144678543$\\
   \hline
\end{tabular}}\label{tab2}
\end{table}
\end{center}

\section{Moduli via potentials}\label{prel}

The finite element method (FEM) is the standard numerical method
for solving elliptic partial differential equations. Since FEM is
an energy minimization method it is eminently suitable for
problems involving Dirichlet energy. In the context of this paper
where the focus is on domains with zero inner angles at the
vertices, the $hp$-FEM variant is the most efficient one
\cite{bs,schwab}. With proper grading of the meshes
even with uniform polynomial order
exponential convergence can be
achieved even in problems with strong corner singularities.


In this section we give a brief overview of the method and our
implementation \cite{hrv,ht}. Of particular importance is the
possibility to estimate the error in the computed quantity of
interest. For quadrilaterals there exists a natural error
estimate, the so-called \textit{reciprocal relation} which is a
necessary but not sufficient condition for convergence. However, if the reciprocal
relation is coupled with \textit{a posteriori} error estimates, we
can trust the results with high confidence \cite{hno}.

\subsection{Modulus of  quadrilateral and the Dirichlet integral}

Let $Q= (Q; z_1, z_2, z_3, z_4)$ be a quadrilateral and the
boundary $\partial Q= \cup_{k=1}^4\partial Q_k$ where all four
boundary arcs are assumed to be non-degenerate. Consider the
following Dirichlet-Neumann problem already introduced in the
introduction:
\begin{equation} \label{eq:dirichlet}
\left\{\begin{matrix}
    \Delta u& =&\ 0,& \text{on}\ &{\ Q,} \\
    u& =&\ 1,&  \text{on}\ &{\partial Q_1,}\\
    u& =&\ 0, & \text{on}\ &{\partial Q_3,}\\
    \partial u/\partial n&  =&\ 0,& \text{on}\ &{\partial Q_2,}\\
    \partial u/\partial n&  =&\ 0,& \text{on}\ &{\partial Q_4.}\\
\end{matrix}\right.
\end{equation}
%

Assume that $u$ is a (unique) harmonic solution of the
Dirichlet-Neumann problem (\ref{eq:dirichlet}). Then the modulus
of $Q$ is defined as
\begin{equation} \label{eq:qmod}
\mbox{\rm Mod}(Q)= \iint_\symD |\nabla  u|^2\,dx \, dy.
\end{equation}

The equality (\ref{eq:qmod}) shows that the modulus of a
quadrilateral is the Dirichlet integral, i.e., the $H^1$-seminorm
of the potential $u$ squared, or, the energy norm squared, a
quantity of interest which is natural in the FEM setting.

\subsection{Mesh refinement and exponential convergence}
\label{sec:refinement}

The idea behind the $p$-version is to associate degrees of freedom
to topological entities of the mesh in contrast to the classical
$h$-version where it is to mesh nodes only. The shape functions
are based on suitable orthogonal polynomials and their supports
reflect the related topological entity, nodes, edges, faces (in
3D), and interior of the elements. The nodal shape functions
induce a partition of unity.

In many problem classes it can be shown that if the mesh is graded
appropriately the method convergences exponentially in some
general norm such as the $H^1$-seminorm. Moreover, due to the
construction of shape functions, it is natural to have large
curved elements in the mesh without significant increase in the
discretization error. Since the number of elements can be kept
relatively low given that additional refinement can always be
added via elementwise polynomial degree, variation in the boundary
can be addressed directly at the level of the boundary
representation in some exact parametric form.

\begin{figure}[h]
    \centering
    \subfloat[Q1: Mesh; $\alpha = \pi/4$.]{\includegraphics[width=0.45\textwidth]{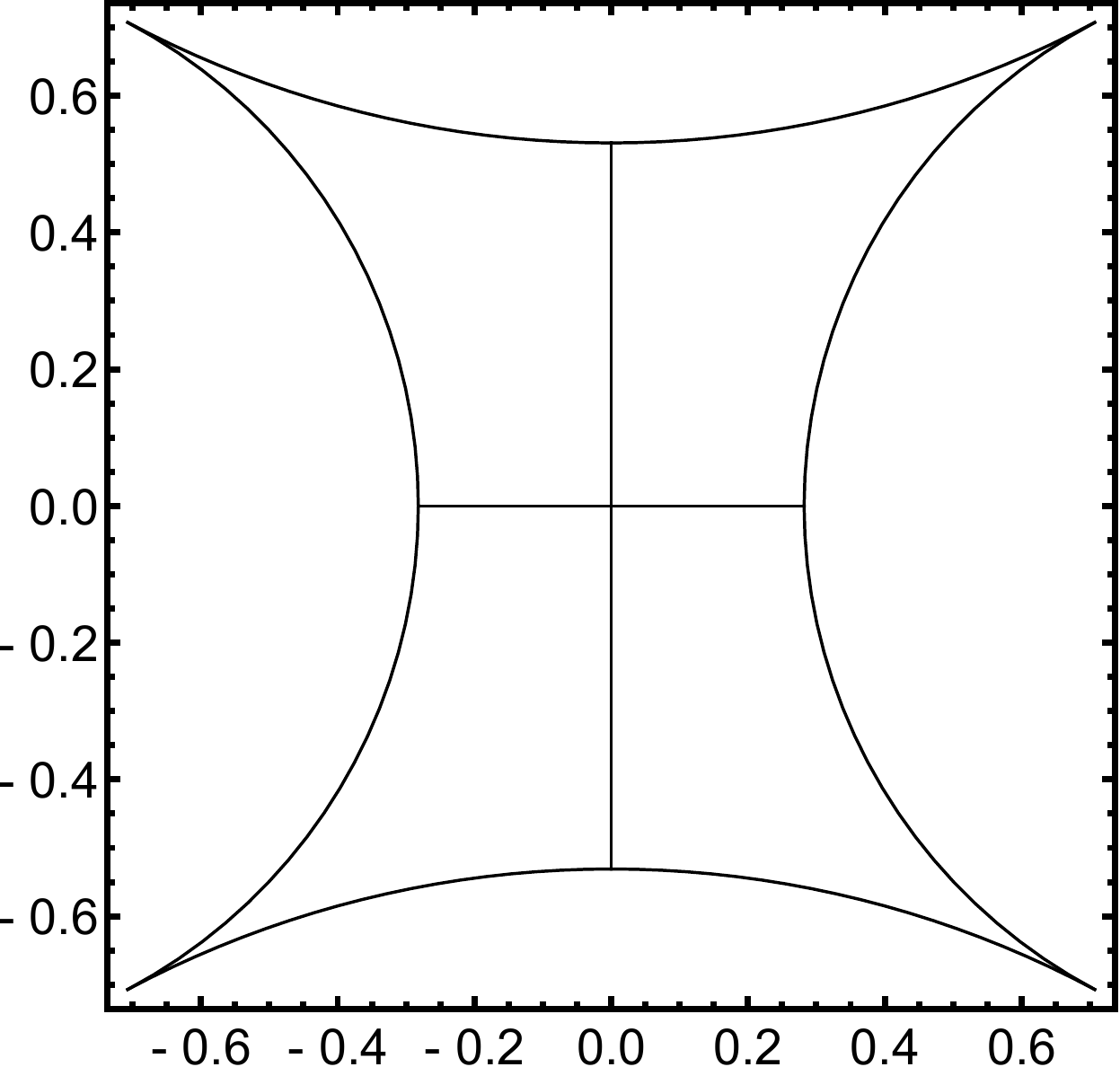}}\quad
    \subfloat[Q2: Mesh; $\alpha = \pi/8$.]{\includegraphics[width=0.45\textwidth]{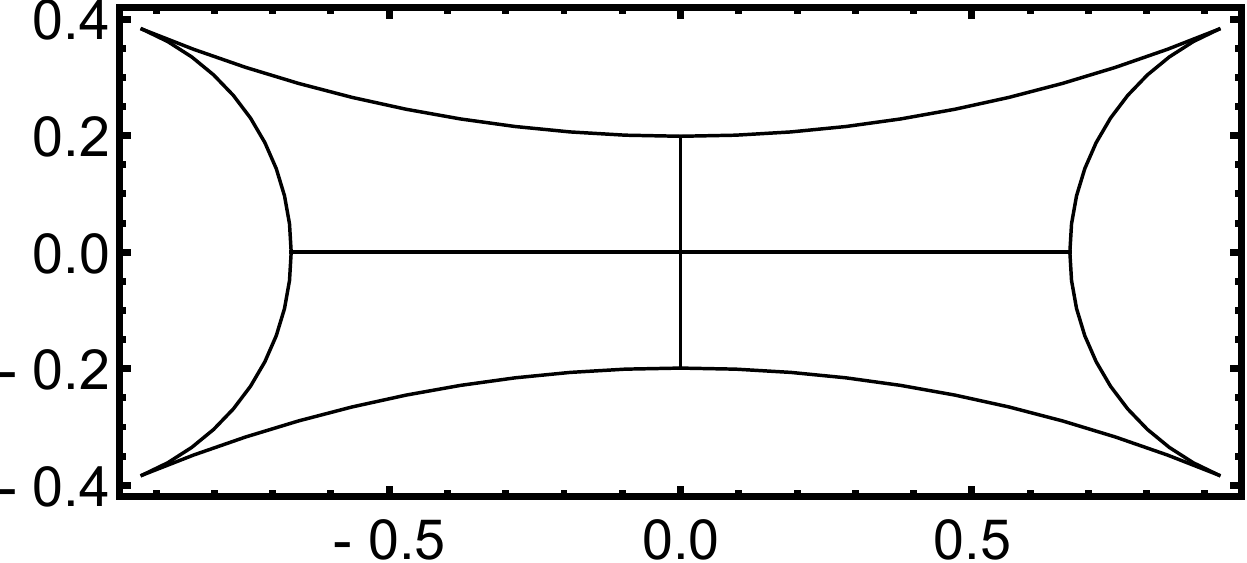}}
    \caption{Circular quadrilaterals. Pure $p$-version meshes.}\label{fig:CQ}
\end{figure}

\begin{figure}[H]
    \centering
    \subfloat[Initial mesh.]{\label{fig:CNGA}\includegraphics[width=0.45\textwidth]{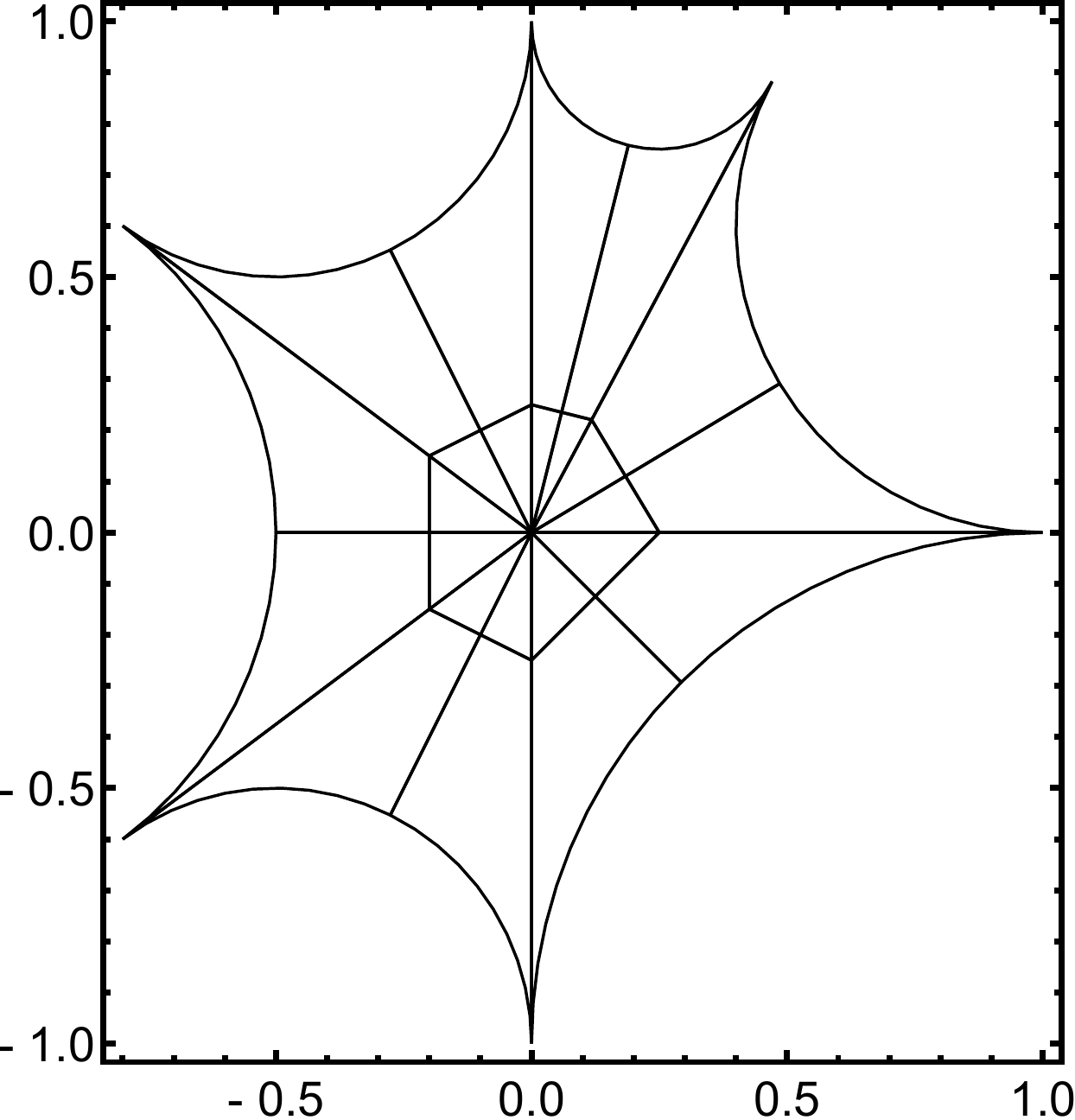}}\quad
    \subfloat[After two levels of refinement.]{\includegraphics[width=0.45\textwidth]{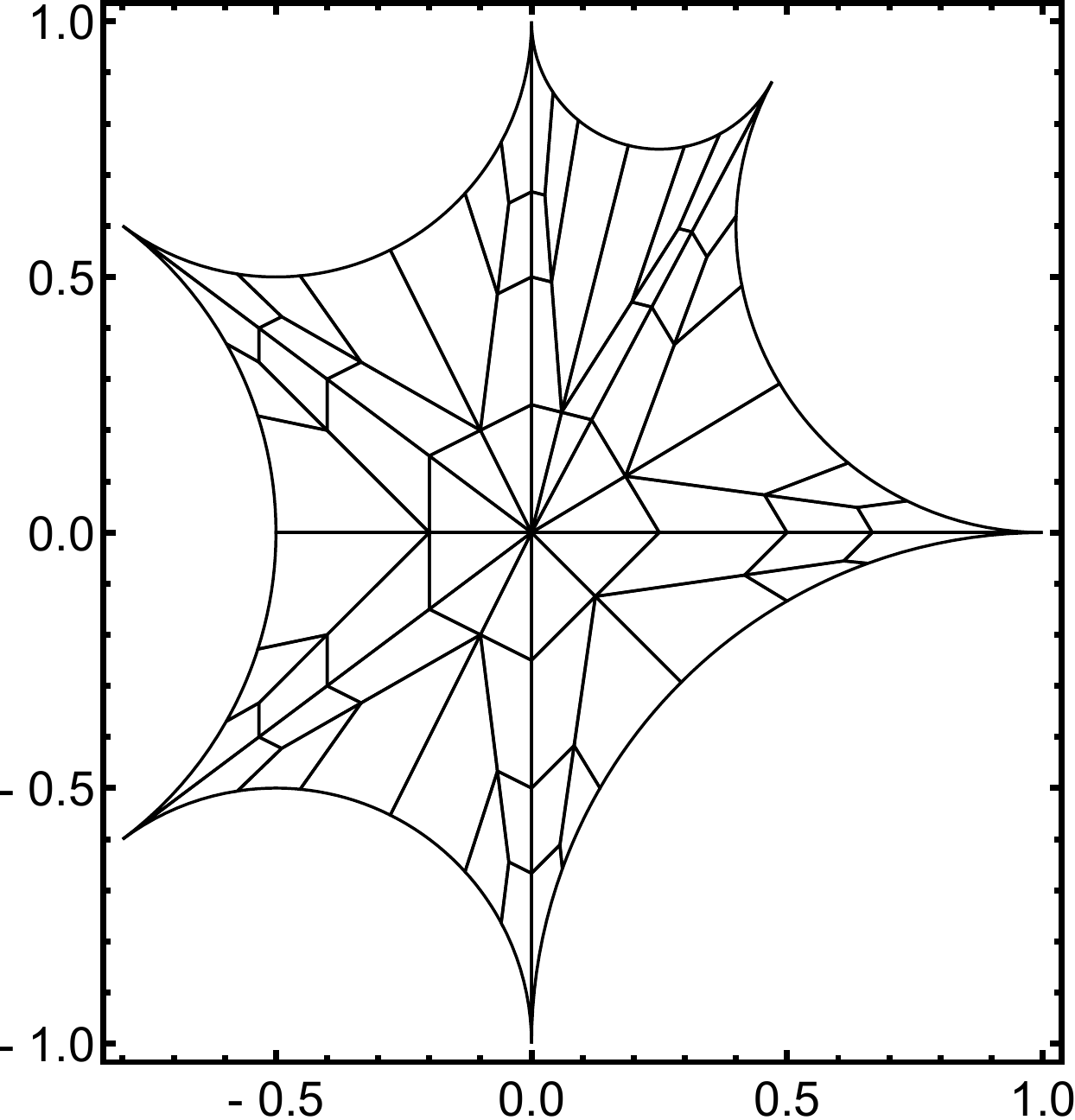}}\\
    \subfloat[Potential of the problem.]{\label{fig:CNGc}\includegraphics[width=0.45\textwidth]{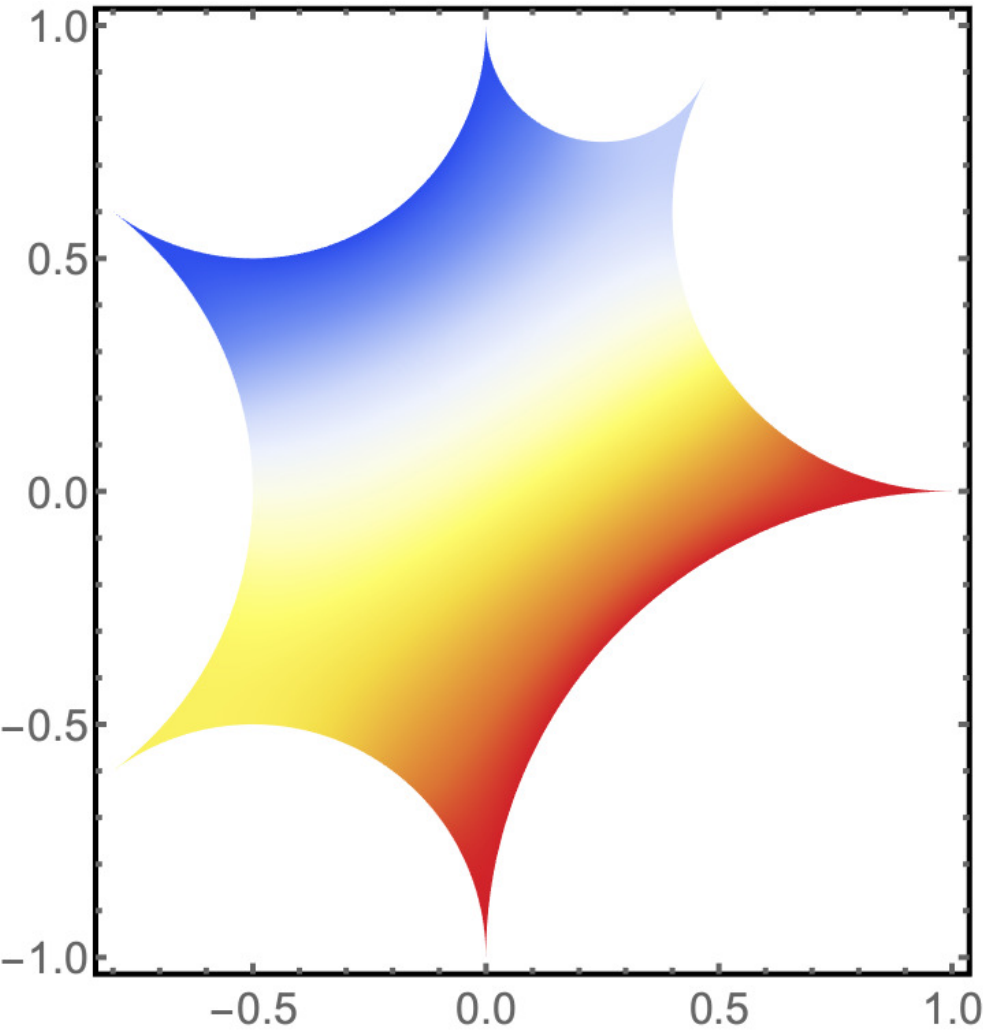}}\quad
    \subfloat[Potential of the conjugate problem.]{\label{fig:CNGd}\includegraphics[width=0.45\textwidth]{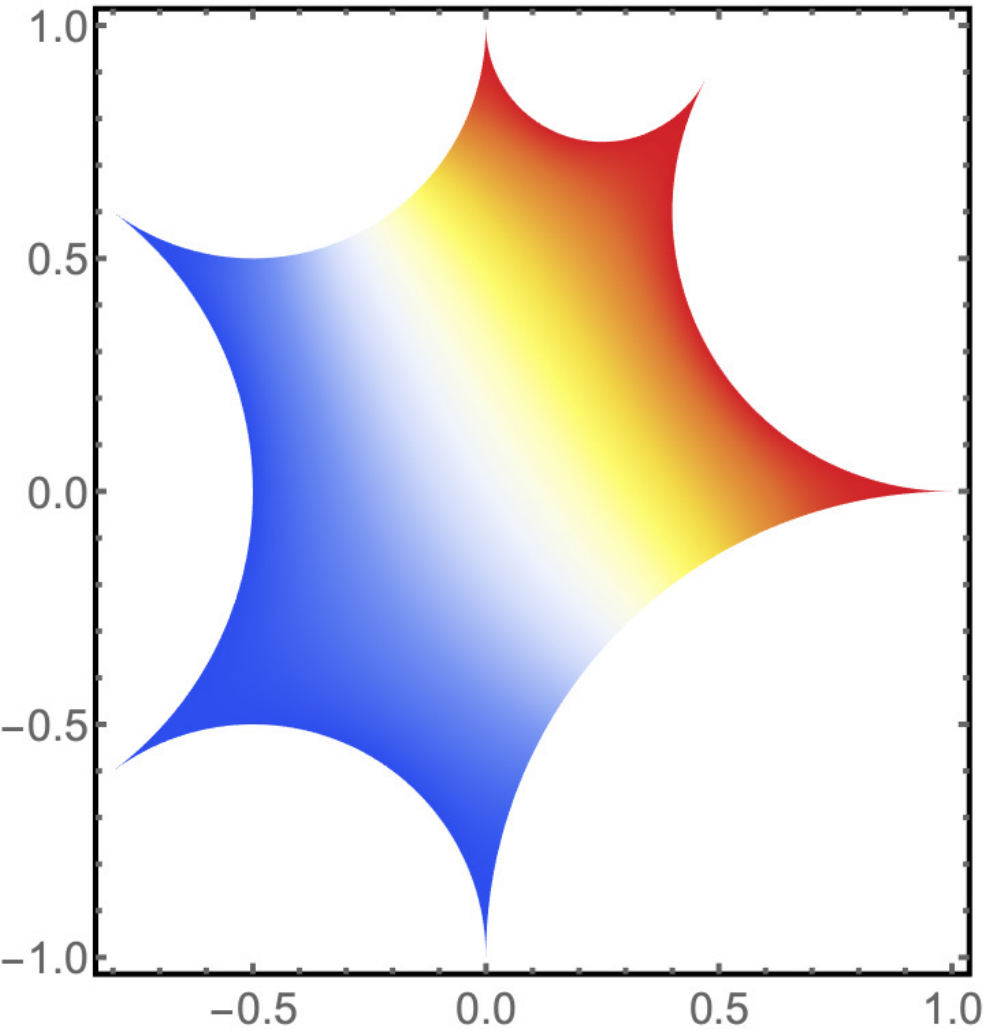}}
    \caption{Circular hexagon. }\label{fig:CNG}
\end{figure}

To fully realize the potential of the $p$-version, one has to
grade the meshes properly and therefore we really use the
$hp$-version here. Consider the meshes in Figures~\ref{fig:CNG}
and \ref{fig:CQ}. In Figures~\ref{fig:CNG} the basic refinement
strategy is illustrated. We start with an initial mesh, where the
corners with singularities are \textit{isolated}, that is, the
subsequent refinements of their neighboring elements do not
interfere with each other. Then the mesh is refined using
successive applications of replacement rules.

In our implementation the geometry can be described in exact
arithmetic and therefore there are not any fixed limits on the
number of refinement levels. In the case of graded meshes one has
to resolve the question of how to set the polynomial degrees at
every element, indeed, a form of refinement of its own. One option
in the case of strong singularities is to set the polynomial
degree based on graph distance from the singularity.
Alternatively, the degree $p$ can be constant over the whole mesh
despite the grading.

\subsection{Error estimation}
\label{sec:errorestimation}

Assuming that the exact capacity is not known, we have two types
of error estimates available: the reciprocal estimate and an a
posteriori estimate. Naturally, if the exact value is known, we
can measure the true error.

\subsubsection{Reciprocal Error Estimate}
The first error estimate is rather unusual in the sense
that it is based on \textit{physics}, yet only necessary. For
every quadrilateral the so-called reciprocal relation can be used
(Definition~\ref{def:recip}). More detailed, from the
definition of modulus via conformal mapping it is clear that the
following \textit{reciprocal identity} holds:
\begin{equation} \label{eq:recip} \mbox{\rm Mod}(\symQuad)\,
\mbox{\rm Mod}(\symQuadC) =1;
\end{equation}
here $\symQuadC= (Q; z_2,$ $z_3,z_4, z_1)$ is called the {\it
quadrilateral} conjugate  to $Q=(Q; z_1,z_2,z_3, z_4)$. In
numerical context (\ref{eq:recip}) gives us the following error
characteristics.

\begin{defn}{\bf (Reciprocal Identity and Error)} \label{def:recip} We will call $$\varepsilon_R = |1-\mbox{\rm Mod}(\symQuad)\,
\mbox{\rm Mod}(\symQuadC)|$$ the error measure and $\varepsilon_N
=\left|\lceil\log_{10}|\varepsilon_R |\rceil\right|$ the related
error number.
\end{defn}

\subsubsection{Auxiliary Space Error Estimate}

Consider the abstract problem setting with $u$ as the standard
piecewise polynomial finite element space on some discretization
$T$ of the computational domain $D$. Assuming that the exact
solution $u \in H_0^1(D)$ has finite energy, we arrive at the
approximation problem: Find $\hat{u} \in V$ such that
\begin{equation}\label{eq:approximation}
  a(\hat{u},v)=l(v)\ (= a(u,v))\quad (\forall v \in V),
\end{equation}where
$a(\,\cdot\,,\,\cdot\,)$ and $l(\,\cdot\,)$, are the bilinear form
and the load potential, respectively. Additional degrees of
freedom can be introduced by enriching the space $V$. This is
accomplished via introduction of an auxiliary subspace or ``error
space'' $W \subset H_0^1(D)$ such that $V \cap W = \{0\}$. We can
then define the error problem: Find $\varepsilon \in W$ such that
\begin{equation}\label{eq:error}
  a(\varepsilon,v)=l(v)- a(\hat{u},v) (= a(u-\hat{u},v))\quad (\forall v \in W).
\end{equation}
This is simply a projection of the residual to the auxiliary space.
In 2D the space $W$, that is, the additional unknowns, can be
associated with element edges and interiors. Thus, for
$hp$-methods this kind of error estimation is natural. The main
result on this kind of estimators is the following
theorem.
\begin{theorem}[\cite{hno}]\label{KeyErrorThm}
There is a constant $K$ depending only on the dimension $d$,
polynomial degree $p$, continuity and coercivity constants $C$ and
$c$, and the shape-regularity of the triangulation $\cT$ such that
\begin{align*}
\frac{c}{C}\,\|\ee\|_1\leq\|u-\hat{u}\|_{1}\leq K
\left(\|\ee\|_{1}+\osc(R,r,\cT)\right),
\end{align*}
where the residual oscillation depends on the volumetric and face
residuals $R$ and $r$, and the triangulation $\cT$.
\end{theorem}

The solution $\varepsilon$ of (\ref{eq:error}) is called the
\textit{error function}. It has many useful properties for both
theoretical and practical considerations. In particular, the error
function can be numerically evaluated and analyzed for any finite
element solution. By construction, the error function is
identically zero at the mesh points. In the examples below, the
space $W$ contains edge shape functions of degree $p+1$ and
internal shape functions of $p+1$ and $p+2$.
This choice is not arbitrary but based on careful cost analysis \cite{hno}.

\subsubsection{Examples}\label{exam}

The numerical examples are defined in Table~\ref{tbl:examples}.
We consider in detail two circular quadrilaterals and one hexagon.
The related results of Table~\ref{tab2} above have been obtained with the
method discussed here.
\begin{table}
  \centering
  \subfloat[Geometry definitions.]{
  \begin{tabular}{ll}
  Example & Parameters or coordinates\\ \hline
  $Q_1$   &  $\vphantom{\int^B}$  $\alpha = \pi/4$, $t=\frac{1}{5} \left(4+\sqrt{2}\right)$, $s=\frac{1}{23} \left(20+19 \sqrt{2}\right)$,
  \\[1mm]
  & $r_1=\frac{1}{5} \sqrt{33-12 \sqrt{2}}$, $r_2=\frac{1}{23} \sqrt{777+300 \sqrt{2}}$  \\[2mm]
  $Q_2$   & $\alpha = \pi/8$, $t=\sec \left(\frac{\pi }{8}\right)$, $s=\csc \left(\frac{\pi }{8}\right)$, \\
  & $r_1=\tan \left(\frac{\pi }{8}\right)$, $r_2=1+\sqrt{2}$ \\[2mm]
    hexagon & $A'=-i$, $B'=1$, $C'=\frac{8}{17}+ \frac{15}{17}\,i$, \\[1mm]
 & $D'=i$, $E'=-\frac{4}{5}+ \frac{3}{5}\,i$, $F'=-\frac{4}{5}- \frac{3}{5}\,i$  \\[2mm] \hline
  \end{tabular}
  }\\
  \subfloat[Problem definitions.]{
  \begin{tabular}{lll}
  Example & Problem & Conjugate\\ \hline
  hexagon & $Q = (Q; A', B', D', E')$ & $\vphantom{\int^B}$ $\symQuadC = (Q; B', D', E', A')$ \\ \hline
  \end{tabular}
  }
  \caption{Examples: Parameters or coordinates used to define the problems.  The problem for hexagon is obtained from Example~\ref{exam1} (cf. Figure~\ref{fig:ngon} (A)) after
carrying out the M\"obius transformation of
Remark~\ref{rem2}.}\label{tbl:examples}
\end{table}

\begin{figure}[H]
    \centering
    \subfloat[{Q1: Reciprocal error; $c = 1/3$.}]{\includegraphics[width=0.4\textwidth]{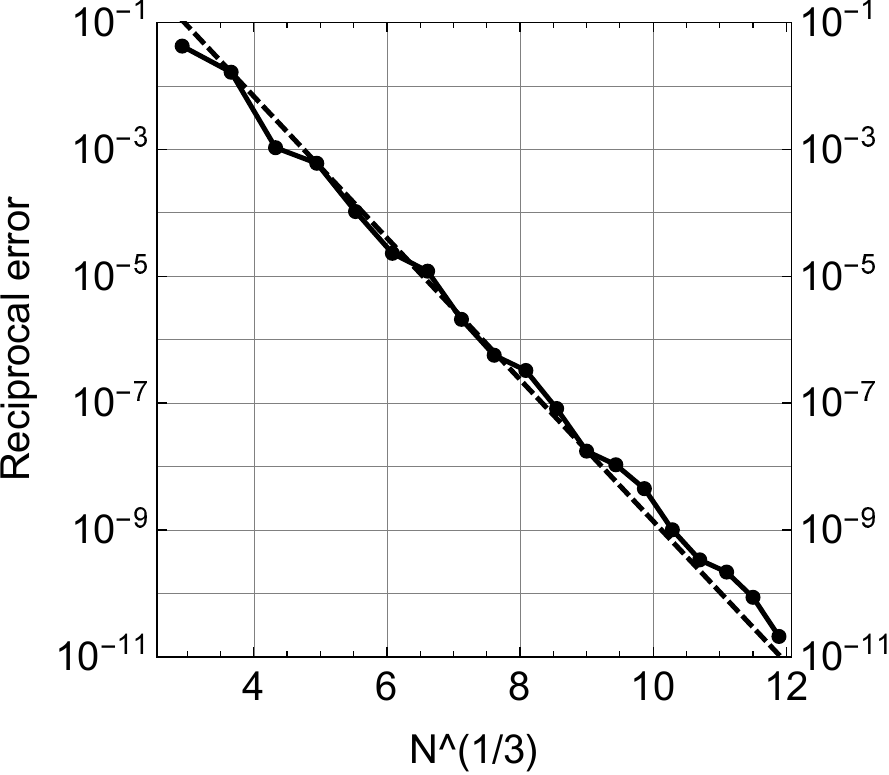}}\quad
    \subfloat[{Q1: Estimated error; $c = 1/3$.}]{\includegraphics[width=0.4\textwidth]{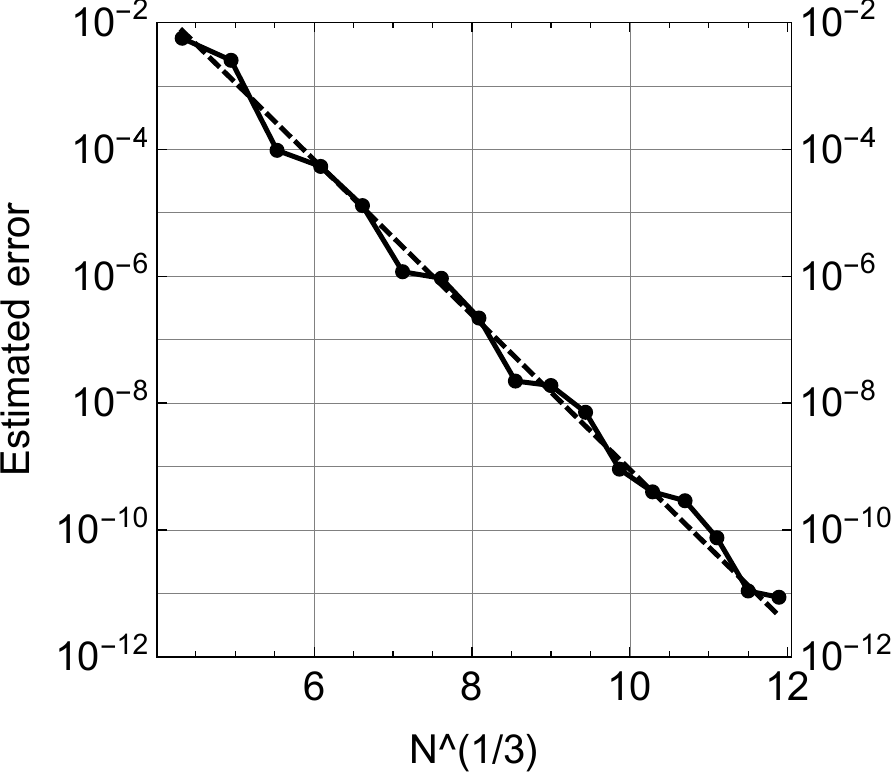}}\\
    \subfloat[{Q2: Reciprocal error; $c = 1/4$.}]{\includegraphics[width=0.4\textwidth]{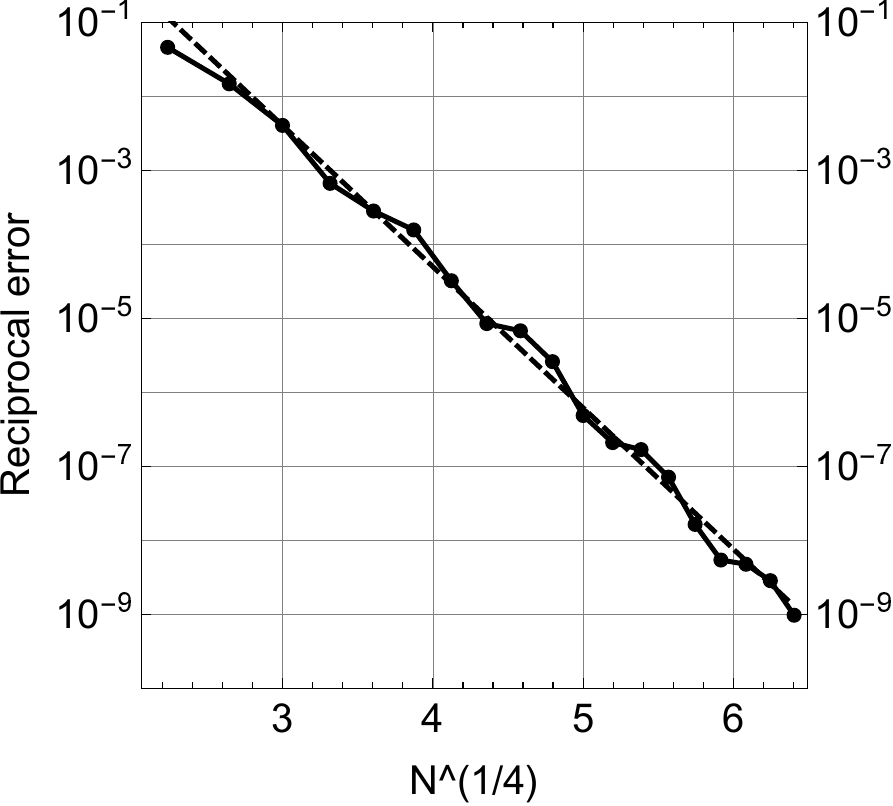}}\quad
    \subfloat[{Q2: Estimated error; $c = 1/4$.}]{\includegraphics[width=0.4\textwidth]{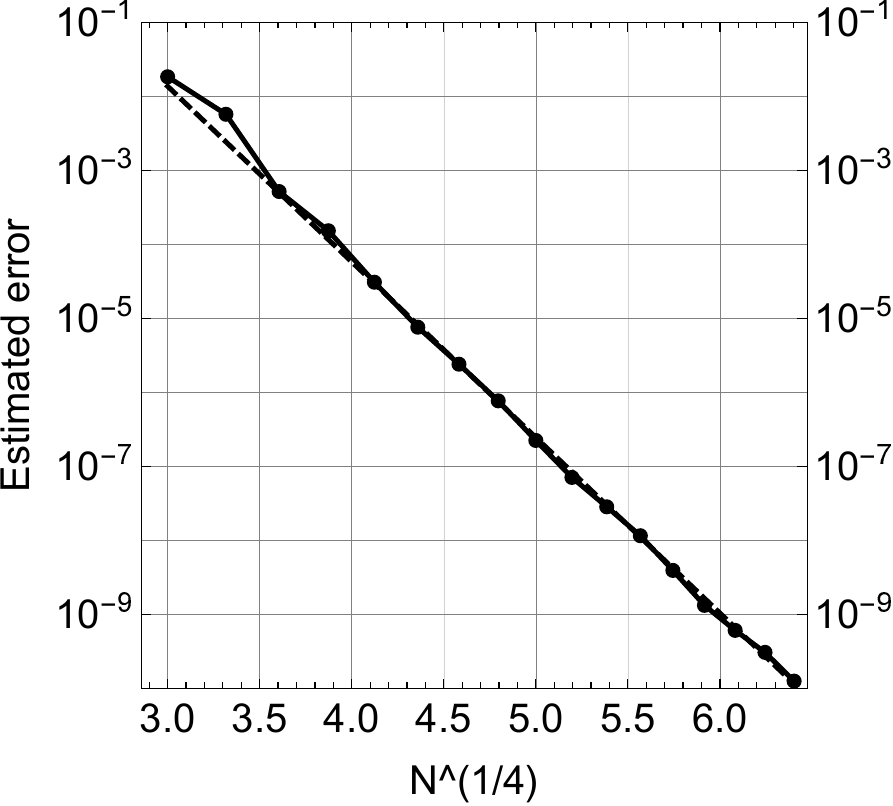}}\\
    \caption{Circular quadrilaterals. Different types of errors vs the number of degrees of freedom (log-plots).
    Solid line with markers represents the observed errors and the dashed line fitted exponential curve (rate $c$ indicated in the caption).}\label{fig:QConvergence}
\end{figure}

\begin{figure}[H]
    \centering
    \subfloat[{Error in capacity; $c = 1/4$.}]{\includegraphics[width=0.4\textwidth]{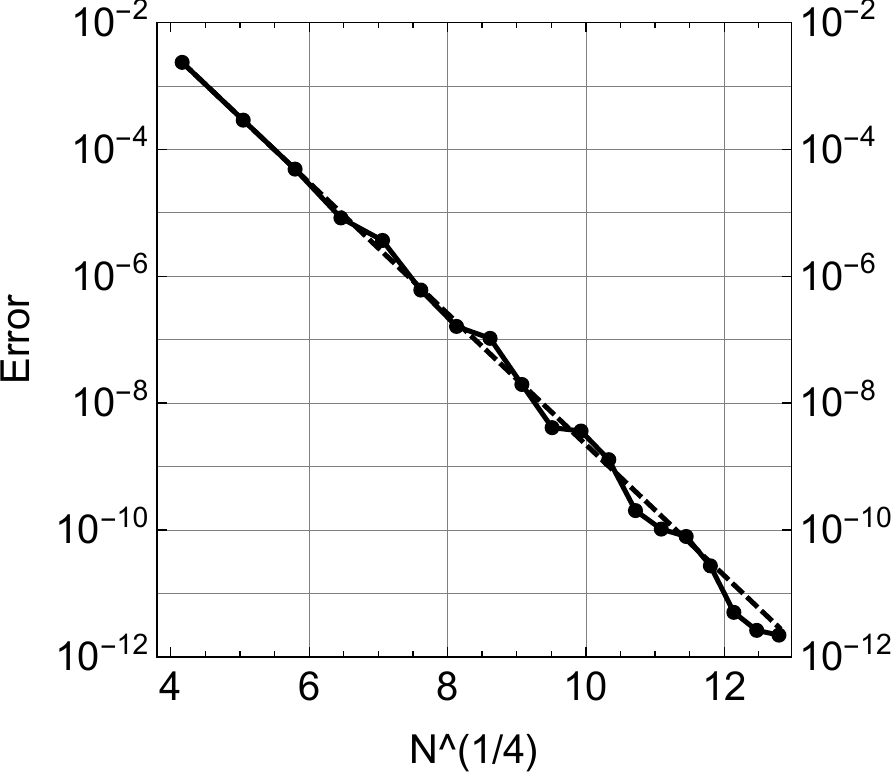}}\quad
    \subfloat[{Reciprocal error; $c = 1/4$.}]{\includegraphics[width=0.4\textwidth]{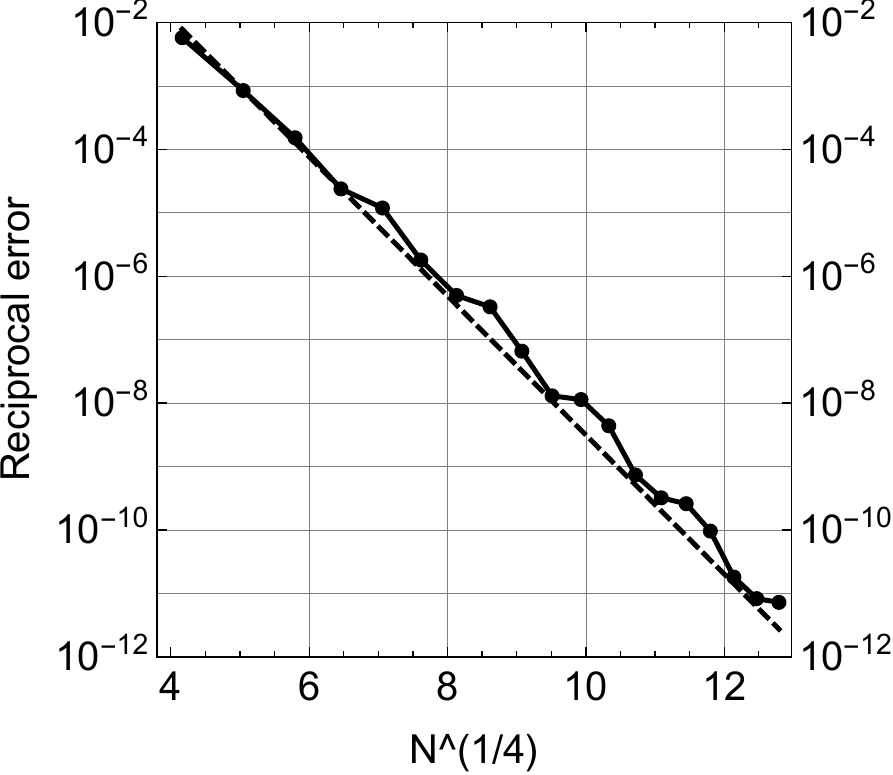}}\\
    \subfloat[{Estimated error; $c = 1/4$.}]{\includegraphics[width=0.4\textwidth]{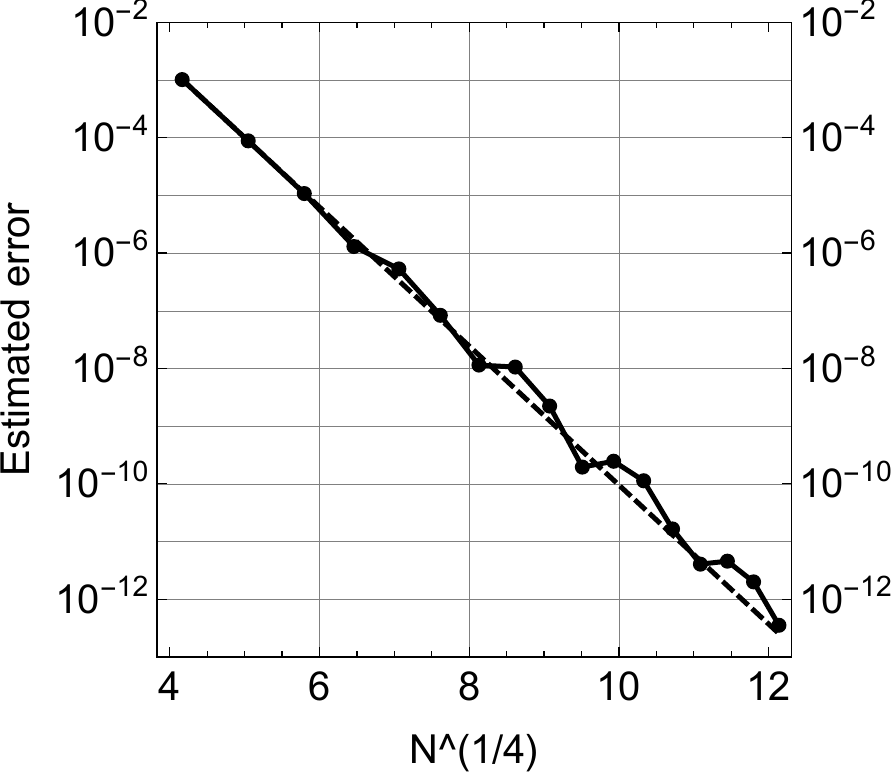}}\quad
    \subfloat[{Comparison of the real and estimated errors.}]{\includegraphics[width=0.4\textwidth]{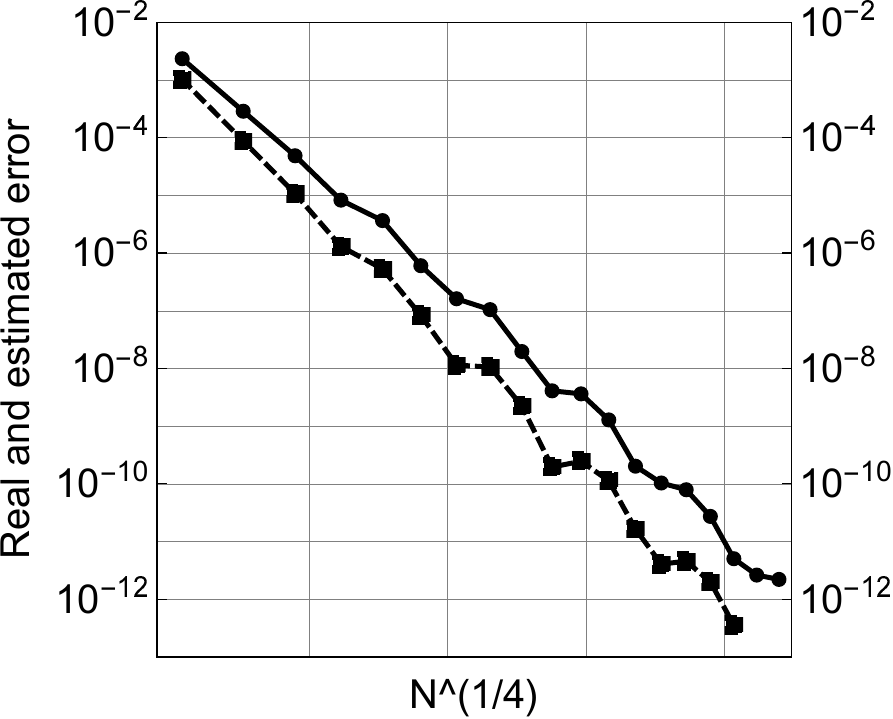}}
    \caption{Circular hexagon. Different types of errors vs the number of degrees of freedom (log-plots).
    Solid line with markers represents the observed errors and the dashed line fitted exponential curve (rate $c$ indicated in the caption).
    In the comparison graph the estimated error is smaller with asymptotic effectivity of 1/10.}\label{fig:CNGConvergence}
\end{figure}

In Figures~\ref{fig:QConvergence} and \ref{fig:CNGConvergence}
different error measures and the related convergence in $p$ are
shown. 
In all cases exponential convergence is realized. The
estimated rates obtained through nonlinear fitting are also
indicated (with dashed lines) and the parameters used are given in
Table~\ref{tbl:fits}. One has to remember that such fits are
notoriously sensitive to selected points and thus, the rates given
here should be taken as possible rates rather than the definitive
ones. We have used the visualization technique where the scaling
is selected to be such that the observed graph appears linear.

For the two quadrilaterals $Q_1$ and $Q_2$ the results are very
good indeed. In fact, the estimated rate for the $Q_1$ is the
theoretically optimal one in terms of the number of degrees of
freedom, of course with a large constant \cite{schwab}. The reason
behind such a spectacular accuracy is that the underlying mapping of
the curved elements, the blending function mapping, is exact for
circular boundary segments. From the point of view of the method
the corner singularity is practically removed by the mapping. For
the $Q_2$ with a smaller aspect ratio, at higher polynomial orders
there is degradation of the convergence rate in comparison to the
symmetric domain of $Q_1$, and indeed, the selected scaling cannot
remain the same as for $Q_1$. Notice that the exceedingly large
constant $a_1$ for the estimated error is due to the non-trivial
oscillation in the estimate and asymptotic convergence reached
only at high values of $p$. Due to symmetry, $\mbox{\rm Mod}(Q_1)
= 1$. We have not shown the error in capacity in this case. It is
done in the following case, however.

For the circular hexagon the geometric meaning of the domain and
its conjugate is illustrated in Figures~\ref{fig:CNGc} and
\ref{fig:CNGd}. In this case the exact modulus is also known,
$$
\mbox{\rm Mod}(Q) = \tau(\sqrt{2})/2 = 
\frac{K(1/\sqrt{1+\sqrt{2}})}{K(\sqrt{\sqrt{2}/(1+\sqrt{2})})}\approx
0.92401502327430725964\ldots,
$$
where $K(r)$ is the complete elliptic integral. The method is indeed very accurate and the
observed rate is within the expected range. As is often the case,
the auxiliary space estimate trails the true error, yet the
effectivity is still over 1/10 as $p$ increases. Due to the
construction of the auxiliary space, the estimate is computed at a
lower polynomial order.

In Table~\ref{tbl:fits} the error numbers for the reciprocal error estimates
are also reported for the highest polynomial order.
These results are aligned with those reported on similar problems before \cite{hrv3}.

\begin{table}
  \centering
  \begin{tabular}{llrrrrr}
  Example & Error Type & $a_1$ & $a_2$ & $c$ & $N$ $(p=20)$ & Error Number \\ \hline
  $Q_1$   & Estimated  & 1397 & -2.8 & 1/3 & 1681 & \\
          & Reciprocal & 198 & -2.6 & 1/3 & 1681 & 10\\
  $Q_2$   & Estimated  & 179191 & -5.5 & 1/4 & 1681 & \\
          & Reciprocal & 2154 & -4.4 & 1/4 & 1681 & 9\\
  hexagon & True       & 47   & -2.4 & 1/4 & 26761 & \\
          & Estimated  & 105  & -2.8 & 1/4 &  $(p=18)$ 21709 & \\
          & Reciprocal & 287   & -2.5 & 1/4 & 26761 & 11\\ \hline
  \end{tabular}
  \caption{Parameters of the nonlinear fits: $a_1 \operatorname{exp}({a_2 N^c})$,
  where $N$ is the number of degrees of freedom, and the error numbers for the reciprocal errors.
  }\label{tbl:fits}
\end{table}

\begin{remark}[On Computational Complexity of the Error Estimates]
The two error estimates do not differ in their computational complexity
in any significant way.
Although the reciprocal error estimate requires the solution of two problems,
and the auxiliary space estimate is for one problem only, the cost of numerical
integration (always an issue in high order methods) is roughly the same, and
the two solution steps for the reciprocal error estimate can share
the Cholesky factorization of the interior degrees of freedom.
\end{remark}

\begin{remark}[On Performance Comparison Between Schwarz ODE and $hp$-FEM]
As mentioned above, the quadrilateral example is particularly well-suited to
the $hp$-FEM. This makes it somewhat awkward to compare the computational
efficiency of the two numerical approaches presented in this paper.
Of course, one should also take into account the time spent in defining the
computational domain. This is very difficult to measure, however.
The Schwarz ODE routine and the $hp$-solver have comparable performance
when the former is run using standard precision. Due to the implementation
of the ODE solver, the higher accuracy is obtained only by changing the
floating point representation which leads to longer run times. On the other hand,
the Schwarz ODE has an almost uniform runtime characteristics over all
circular quadrilaterals and it is likely that replacing the general ODE solver routines with problem specific
ones will lead to significant improvements in run times.
The $hp$-FEM requires more resources if the
discretization includes more elements. With the current $hp$-implementation the
non-graded discretization of the $n$-gon (Figure~\ref{fig:CNGA}) took three times longer than the
corresponding quadrilaterals (Figure~\ref{fig:CQ}).
\end{remark}
\section{Conclusions}
\label{sec:conclusions} Here moduli of planar circular
quadrilaterals symmetric with respect to both the coordinate axes
have been investigated. Computation of moduli of planar domains
with cusps is difficult and requires either a customized, analytic
algorithm or general method with sufficient flexibility. The
Schwarz ODE introduced here, an analytic method to determine a
conformal mapping the unit disk onto a given circular
quadrilateral, belongs to the first category. $hp$-FEM on the
other hand provides a framework for highly efficient numerical PDE
solvers. We have shown that these two different approaches provide
results agreeing with high accuracy over two sets of parametrized
examples.

\appendix
\section{Reference Implementations}\label{impl}
The programs used to compute the Table~\ref{tab2} are available at
\begin{center}
 \texttt{https://github.com/hhakula/hnv}
\end{center}
and
  Version 1.0, used in this paper, is archived at DOI: \texttt{10.5281/zenodo.4718320}.

 The Schwarz ODE code is also listed below. The expected output of
the program is $$\{1.02791, 0.440765, 1.25503\}.$$
\begin{lstlisting}[language=Mathematica,caption={Schwarz ODE}]
   t = 2.0174131664886366`;
   s = 1.1416407864998739`;
   r1 = 1.642663833605752`;
   r2 = 0.6753740370343625`;
   Clear[K, b, gamma];
   K = r2/r1;
   k = s/t;
mod = 1.;
While[Abs[mod - 1.] < 10^(-5),
ba = 0.;
bb = Pi/4;
  Do[Clear[gamma,b];
     b = (ba + bb)/2.;
     theta1 = 0;
     theta2 = b/2.;
     theta3 = Pi/2*(1/2) + b*(1/2);
     theta4 = Pi/2.;
     agamma = 0.7 - 4./Pi*b;
     bgamma = 1.2 - 3./Pi*b;
     F[x_, beta_, theta_, gamma_] = Exp[2*I*theta](Exp[2*I*beta]/(x^2*Exp[2*I*theta]
     -Exp[2*I*beta])^2 +Exp[-2*I*beta]/(x^2*Exp[2*I*theta] - Exp[-2*I*beta])^2
     - gamma (1/((x^2*Exp[2*I*theta] - Exp[2*I*beta])(x^2*Exp[2*I*theta]
     -Exp[-2*I*beta]))));
     Do[gamma = (agamma + bgamma)/2.;
        sol3 = NDSolve[{
         u3''[x]+F[x,b,theta3,gamma]*u3[x]==0,v3''[x]+F[x,b,theta3,gamma]*v3[x]==0,
         u3[0]==0,u3'[0]==Exp[I*theta3],v3[0]==1,v3'[0]==0},{u3,v3},{x,0,1}];
        sol4 = NDSolve[{
         u4''[x]+F[x,b,theta4,gamma]*u4[x]==0,v4''[x]+F[x,b,theta4,gamma]*v4[x]==0,
         u4[0]==0,u4'[0]==Exp[I*theta4],v4[0]==1,v4'[0]==0},{u4,v4},{x,0,1}];
        F3[x_] = u3[x]/v3[x] /. sol3;
        F4[x_] = u4[x]/v4[x] /. sol4;
        x2 = Re[F3[1]];
        y2 = Im[F3[1]];
        y3 = Im[F4[1]];
        S = (1/2)(y2+y3+x2^2/(y2-y3));
        If[ S[[1]]<0,
            bgamma = gamma,
            agamma = gamma
        ],{i, 10}];
     BGAMMA = gamma;
     Clear[gamma];
     agamma = 0.75-4./Pi*b;
     bgamma = 1.2-3./Pi*b;
     Do[gamma = (agamma + bgamma)/2.;
        sol1 = NDSolve[{
        u1''[x]+F[x,b,theta1,gamma]*u1[x]==0,v1''[x]+F[x,b,theta1,gamma]*v1[x]==0,
        u1[0]==0,u1'[0]==Exp[I*theta1],v1[0]==1,v1'[0]==0},{u1, v1},{x,0,1}];
        sol2 = NDSolve[{
        u2''[x]+F[x,b,theta2,gamma]*u2[x]==0,v2''[x]+F[x,b,theta2,gamma]*v2[x]==0,
        u2[0]==0,u2'[0]==Exp[I*theta2],v2[0]==1,v2'[0]==0},{u2,v2},{x,0,1}];
        F1[x_] = u1[x]/v1[x] /. sol1;
        F2[x_] = u2[x]/v2[x] /. sol2;
        x1 = Re[F1[1]];
        x2 = Re[F2[1]];
        y2 = Im[F2[1]];
        T = (1/2)(x1+x2+y2^2/(x2-y1));
        If[ T[[1]] > 0,
            bgamma = gamma,
            agamma = gamma
        ], {i, 10}];
     AGAMMA = gamma;
     Clear[gamma,x1,x2,y2,x3,y3,y4];
     Do[gamma = (AGAMMA + BGAMMA)/2.;
        sol1 = NDSolve[{
        u1''[x]+F[x,b,theta1,gamma]*u1[x]==0,v1''[x]+F[x,b,theta1,gamma]*v1[x]==0,
        u1[0]==0,u1'[0]==Exp[I*theta1],v1[0]==1,v1'[0]==0},{u1,v1},{x,0,1},
        PrecisionGoal->15];
        sol2 = NDSolve[{
        u2''[x]+F[x,b,theta2,gamma]*u2[x]==0,v2''[x]+F[x,b,theta2,gamma]*v2[x]==0,
        u2[0]==0,u2'[0]==Exp[I*theta2],v2[0]==1,v2'[0]==0},{u2,v2},{x,0,1},
        PrecisionGoal->15];
        sol3 = NDSolve[{
        u3''[x]+F[x,b,theta3,gamma]*u3[x]==0,v3''[x]+F[x,b,theta3,gamma]*v3[x]==0,
        u3[0]==0,u3'[0]==Exp[I*theta3],v3[0]==1,v3'[0]==0},{u3,v3},{x,0,1},
        PrecisionGoal->15];
        sol4 = NDSolve[{
        u4''[x]+F[x,b,theta4,gamma]*u4[x]==0,v4''[x]+F[x,b,theta4,gamma]*v4[x]==0,
        u4[0]==0,u4'[0]==Exp[I*theta4],v4[0]==1,v4'[0]==0},{u4,v4},{x,0,1},
        PrecisionGoal->15];
        F1[x_] = u1[x]/v1[x] /. sol1;
        F2[x_] = u2[x]/v2[x] /. sol2;
        F3[x_] = u3[x]/v3[x] /. sol3;
        F4[x_] = u4[x]/v4[x] /. sol4;
        x1 = Re[F1[1]];
        x2 = Re[F2[1]];
        y2 = Im[F2[1]];
        T = (1/2)(x1+x2+y2^2/(x2-x1));
        R1 = Abs[x1-T];
        x3 = Re[F3[1]];
        y3 = Im[F3[1]];
        y4 = Im[F4[1]];
        S = (1/2)(y3+y4+x3^2/(y3-y4));
        R2 = Abs[y4-S];
        If[ k*T[[1]] < S[[1]],
            BGAMMA = gamma,
            AGAMMA = gamma
        ],{i,25}];
     If[ R2[[1]]/R1[[1]]<K,
         bb = b,
         ba = b
     ],{n,25}];
       m = (Tan[b/2])^4;
       mod = 2*EllipticK[m]/EllipticK[1-m];
       If[ Abs[mod-1.]<10^(-5),
           K = r1/r2
       ];
       If[ Abs[mod-1.]<10^(-5),
           k = t/s
       ];
       ]; (* End top While *)
       If[ K==r2/r1,
           b1 = b,
           b1 = Pi/2 - b
       ];
      m1 = (Tan[b1/2])^4;
      modQ = 2*EllipticK[m1]/EllipticK[1-m1];
      {b1, gamma, modQ}
  \end{lstlisting}

\end{document}